\documentclass{pspum-l}

\usepackage{amssymb}
\usepackage{latexsym}
\usepackage{amsxtra}
\usepackage{amscd}
\usepackage{epic}
\usepackage{eepic}

\makeatother

\newtheorem{thm}{Theorem}[section]
\newtheorem{lem}[thm]{Lemma}
\newtheorem{prop}[thm]{Proposition}
\newtheorem{cor}[thm]{Corollary}

\theoremstyle{definition}
\newtheorem{defn}[thm]{Definition}
\newtheorem{exa}[thm]{Example}
\newtheorem{cond}[thm]{Condition}

\theoremstyle{remark}
\newtheorem{rem}[thm]{Remark}

\numberwithin{equation}{section}

\newcommand{\ZZ}{\mathbb{Z}}
\newcommand{\QQ}{\mathbb{Q}}
\newcommand{\CC}{\mathbb{C}}
\newcommand{\RR}{\mathbb{R}}
\newcommand{\PP}{\mathbb{P}}

\newcommand{\FF}{\mathbb{F}}

\newcommand{\ZZh}{\hat{\ZZ}}

\newcommand{\U}{\mathcal{U}}

\renewcommand{\O}{\mathcal{O}}
\newcommand{\Od}{\hat{\O}}
\newcommand{\Hu}{\mathcal{H}}
\renewcommand{\o}{\mathfrak{o}}

\newcommand{\Xb}{\bar{X}}

\newcommand{\Dt}{\tilde{D}}

\renewcommand{\tt}{\tilde{t}}

\newcommand{\zt}{\tilde{z}}
\newcommand{\zb}{\bar{z}}

\newcommand{\At}{\tilde{A}}

\newcommand{\Spec}{\mathop{\rm Spec}} 
\newcommand{\Hom}{\mathop{\rm Hom}\nolimits}

\newcommand{\Gal}{\mathop{\rm Gal}}

\newcommand{\Inn}{\mathop{\rm Inn}} 
\renewcommand{\Im}{\mathop{\rm Im}} 
\renewcommand{\Re}{\mathop{\rm Re}}

\newcommand{\ord}{\mathop{\rm ord}}
\newcommand{\Ni}{\mathop{\rm Ni}\nolimits}
\newcommand{\cha}{\mathop{\rm char}}
\newcommand{\Br}{\mathop{\rm Br}\nolimits}
\newcommand{\pre}{\mathop{\rm pre}}
\newcommand{\Frac}{\mathop{\rm Frac}}

\renewcommand{\phi}{\varphi}
\newcommand{\phib}{\bar{\phi}}

\newcommand{\psib}{\bar{\psi}}

\newcommand{\nr}{^{\rm\scriptscriptstyle nr}}
\newcommand{\inn}{^{\rm\scriptscriptstyle in}}

\newcommand{\et}{_{\rm\scriptscriptstyle \acute{e}t}}

\newcommand{\inj}{\hookrightarrow}
\newcommand{\To}{\;\longrightarrow\;}
\newcommand{\iso}{\stackrel{\sim}{\to}}
\newcommand{\liso}{\;\stackrel{\sim}{\longrightarrow}\;}

\newcommand{\lpfeil}[1]{\stackrel{#1}{\To}}

\newcommand{\op}[2]{\sideset{^{#1}}{}{\mathop{#2}}}

\newcommand{\g}{{\bf g}}

\newcommand{\mmu}{\boldsymbol{\mu}}
\newcommand{\p}{\mathfrak{p}}

\newcommand{\QQb}{\bar{\QQ}}
\newcommand{\Gb}{\bar{G}}

\newcommand{\fb}{\bar{f}}
\newcommand{\hb}{\bar{h}}
\newcommand{\gb}{\bar{g}}

\newcommand{\gen}[1]{\mathopen<#1\mathclose>}

\newcommand{\rupfeil}[1]{\Big\downarrow\vcenter{%
                         \rlap{$\scriptstyle #1$}}}

\begin{document}

\title{Field of moduli and field of definition of Galois covers}

\author{Stefan Wewers}
\address{University of Pennsylvania}
\email{wewers@math.upenn.edu}


\date{June 2000}

\begin{abstract}
  In this paper we investigate the cohomological obstruction for the
  field of moduli of a $G$-cover to be a field of definition, in the
  case of local fields and covers with tame admissible reduction.
  This applies in particular to $p$-adic fields where $p$ does not
  divide the order of the group $G$.  We give examples of $G$-covers
  with field of moduli $\QQ_p$ that cannot be defined over $\QQ_p$,
  for all primes $p>5$.
\end{abstract}

\maketitle


\section*{Introduction}

In the context of the regular inverse Galois problem, the method of
{\em rigidity} and its generalizations using the braid action and
Hurwitz spaces have been very successful. One drawback of these
methods is that they apply a priori only to groups with trivial
center. To give an example, let $G$ be a finite group and
$f:Y\to\PP^1$ a $G$-Galois cover of the projective line defined over
$\QQb$. Assume that $f$ has rational branch points
$z_1,\ldots,z_r\in\PP^1(\QQ)$ and that the associated tuple of
conjugacy classes $(C_1,\ldots,C_r)$ is rational and rigid. Then $\QQ$
is the {\em field of moduli} of the $G$-cover, i.e.\ for each
$\sigma\in\Gal(\QQb/\QQ)$ the conjugate cover $\op{\sigma}{f}$ is
isomorphic to $f$. If, moreover, $G$ has trivial center, then $\QQ$ is a
{\em field of definition} of $f$, i.e.\ there exists a model
$f_{\QQ}:Y_{\QQ}\to\PP^1_{\QQ}$ of $f$ over $\QQ$. This yields the
basic Rigidity Criterion, see e.g.\ \cite{MM}, \cite{SerreTopics} or
\cite{Voelklein}. 

However, if the center $C$ of $G$ is not trivial, the field of moduli
$k_m$ of a $G$-cover need not be a field of definition. There is an
obstruction $\omega\in H^2(k_m,C)$ such that $f$ can be defined over
an extension $k/k_m$ if and only if the restriction of $\omega$ to $k$
vanishes. This cohomological approach has been developed in several
papers, in the context of Galois covers and in more general
situations, see e.g.\ \cite{CooHar85}, \cite{DebDou1}. In
\cite{DebDou2} it is used to prove a local-to-global principle: a
$G$-cover with field of moduli $\QQ$ is defined over $\QQ$ if and only
if it can be defined over $\QQ_p$, for each prime number $p$. In
\cite{DebHar98} it is shown that $\omega|_{\QQ_p}$ vanishes for each
prime number $p$ at which the cover $f$ has good reduction; this
includes all ``good'' primes in the sense of \cite{Beckmann89}. In the
present paper we describe a method that allows, in many cases, to
explicitly compute the local obstruction $\omega|_{\QQ_p}$, at a prime
$p$ at which the cover $f$ has bad reduction (but which does not
divide the order of $|G|$). We also apply this method in a nontrivial
example.

\vspace{2ex} To explain our results in more detail, let $k$ be a field
which is henselian with respect to a discrete valuation $v$, and let
$f:Y\to\PP^1$ be a $G$-cover, defined over the algebraic closure of
$k$, with field of moduli $k$. Our main assumption is that $f$ has
{\em tame admissible} reduction. For instance, this condition holds if
the order of $G$ is prime to the residue characteristic of $k$. In
this case, bad reduction is caused by degeneration of the branch locus
of the cover modulo $v$. In particular, the results we present here
are interesting only for covers with at least $4$ branch points. By
standard lifting results, the special fiber of a tame admissible
cover, together with certain degeneration data, `knows everything'
about the generic fiber, see e.g.\ \cite{deform}.  Therefore, one can
hope to compute the obstruction for $k$ to be a field of definition
purely in terms of the special fiber of $f$.

Consider the following short exact sequence of cohomology groups:
\begin{equation} \label{introeq1}
  0 \to H^2(k_0,C) \To H^2(k,C) \lpfeil{r_v} H^1(k_0,C(-1)) \to 0,
\end{equation}
where $k_0$ denotes the residue field of $k$, see \cite{SerreCG},
Section II.A.2. Suppose $f:Y\to\PP^1$ is a $G$-cover with field of
moduli $k$ and $\omega\in H^2(k,C)$ is the associated obstruction for
$k$ to be a field of definition. The class $r_v(\omega)\in
H^1(k_0,C(-1))$ is called the {\em residue} of $\omega$ at $v$. If
$r_v(\omega)=0$, then $\omega$, regarded as a class in $H^2(k_0,C)$,
represents the obstruction for the special fiber $\bar{f}$ of $f$ to
be defined over $k_0$ (as an admissible $G$-cover). This happens, for
instance, if $f$ has good reduction. On the other hand, if
$H^2(k_0,C)=0$ (for instance, if $k_0$ is a finite field) then the
special fiber $\fb$ has a $k_0$-model $\fb_{k_0}$, and then
$r_v(\omega)$ represents the obstruction for $\fb_{k_0}$ to lift to a
$k$-model of $f$. These considerations immediately give a new proof of
the main result of \cite{DebHar98}.

If the residue field $k_0$ has characteristic $0$, one can describe
the degeneration behavior of $f$ purely in terms of the {\em Hurwitz
data} attached to $f$. A substantial part of the present paper is
devoted to showing that, in many cases, this description is sufficient
for computing $r_v(\omega)$.

To compute $\omega$ over a $p$-adic field, we use a specialization
technique. Given a $p$-adic field $k$ and a $G$-cover $f$ with field
of moduli $k$, one can construct (under certain conditions) a
$G$-cover $f_t$, with field of moduli $K_t:=k((t))$, such that $f=f_a$
is a `specialization' of $f_t$, for some value $a\in k$ with $v(a)>0$.
Since $k$ has characteristic $0$, one can apply the more geometric
methods referred to in the preceding paragraph to compute the
obstruction $\omega_t$ for $K_t$ to be a field of definition of
$f_t$. In this situation, we prove that the obstruction $\omega$
corresponding to $f$ can be computed by `specializing' $\omega_t$ to
$t=a$. The proof uses the theory of Hurwitz spaces in mixed
characteristic, see e.g.\ \cite{diss}.

Let $\At_5$ be the unique nonsplit central extension of $A_5$ by $\pm
1$. We construct, for each prime number $p>5$, an $\At_5$ cover with
field of moduli $\QQ_p$ which is not defined over $\QQ_p$. To my
knowledge, this kind of example is essentially new (\cite{CooHar85},
Example 2.6, for instance, uses only the real numbers).

We should compare this paper to several approaches in the literature to
understand the absolute Galois group of $\QQ$ acting on the fundamental
groups of moduli spaces.  
\cite{IharMat}  provides a Hurwitz space context for the
Drinfeld-Ihara-Grothendieck relations (that apply  to elements of the
absolute Galois group; call these DIG relations). The approach was
through {\sl tangential base points\/} . We follow that tradition,
though we are not taking the exact same tangential base points. For
example, we often use complex conjugate pairs of branch points, while
they always used sets of real branch points. Ihara's use of the DIG
relations has been primarily to describe the Lie algebra of the absolute
Galois group acting through various pronilpotent braid groups,
especially on the 3 punctured $\lambda$-line.  Fried in
\cite[App.~C]{FrMT}  proposed Modular Towers, a profinite construction,
as suitably like finite representations of the fundamental group to see
the DIG relations at a {\sl finite level\/} . There is an analogy with
Ihara in that modular curves are close to considerations about the
$\lambda$-line, though there is no phenomenon related to them that
suggests seeing the DIG relations. Still, modular curves are just one
case of Modular Towers.

Fried specifically proposed the Modular Tower attached to $A_5$ and four
3-cycles as a candidate where one would see a system of Serre
obstruction situations from covers of $A_5$.
\cite[Prop.~8.12]{BaileyFriedTV}  uses a specific case of Modular Tower
levels having  a  tower of these Serre obstructions. Over the reals,
there are points, {\sl Harbater-Mumford (HM) reps.\/} ~ at every level
of this Modular Tower where (in our notation) the Serre obstruction
vanishes and has value 0, and other points ({\sl near HM reps.}) where
it vanishes, but does not have value 0. At levels 1 and beyond, these
are the only real points and the real components through them end at
cusps of reduced Hurwitz spaces (covering the $j$-line). These
computations are slightly differently than ours, and they give their own
approach to Serre's obstruction through direct use of rational points on
covers \cite[Prop.~6.8]{BaileyFriedTV} . Still, restricting our
Prop.~\ref{computeprop}  to the reals is essentially equivalent to an
example of theirs.

This suggests there are analogs of $p$-adic points on these Modular
Towers levels that have a similar tangential base point (cusp geometry)
analysis to the near HM and HM reps. Further, this geometry should
reveal the DIG relations on  actual covers, instead of as a Lie algebra
relation. While there may be technical difficulties with carrying this
out at all levels, the computations of \cite[Prop.~9.8]{BaileyFriedTV}
at level 1 of this $A_5$ Modular Tower should be feasible.

\vspace{2ex}
The paper is organized in two sections. In Section 1, we recall
general results about the field of moduli and the field of definition
of Galois covers, from our point of view. Section 2 is concerned with
the case of a henselian ground field and contains the main results.

\vspace{2ex} The author would like to thank the organizers of the
special semester {\em Galois Groups and Fundamental Groups} for
inviting him and for providing financial support. He would also like
to thank the MSRI for its hospitality, and the referee for useful
comments. The final version of this paper was written while the author
received a grant from the {\em Deutsche Forschungsgemeinschaft}.



\section{General results}

In Section \ref{fod} we recall the necessary definitions. In Section
\ref{hurw} we introduce some notation concerning the Hurwitz
description of Galois covers and compute the f.o.d.-obstruction in an
example (the group is ${\rm SL}_2(\ell)$ and the ramification type is
$(4,\ell,\ell)$). In Section \ref{real} we study in some detail Galois
covers over the real numbers.

\subsection{The field of moduli condition and the f.o.d.-obstruction} 
\label{fod}

Let $G$ be a finite group and $k$ a field. We denote by $k^s$ a fixed
separable closure of $k$ and by $\Gamma_k:=\Gal(k^s/k)$ the absolute
Galois group of $k$. Suppose we are given a $G$-cover
\[
       f:Y\lpfeil{G} \PP^1_{k^s}
\]
of the projective line over $k^s$; by this we mean that $f$ is a
finite Galois cover of smooth projective curves over $k^s$, with
Galois group $G$. We will always assume that $f$ is {\em tamely
  ramified}, and we let $S=\{z_1,\ldots,z_r\}\subset\PP^1_{k^s}$ be
the set of branch points of $f$. Thus, the $G$-cover $f$ corresponds
to a surjective morphism $\Phi:\pi_1^t(U)\to G$, where $\pi_1^t(U)$
denotes the tame fundamental group of $U:=\PP^1_{k^s}-S$.

The question we are concerned with is the following: is $k$ a {\em
  field of definition} of $f$, i.e.\ does there exist a $G$-cover
$f_k:Y_k\to\PP^1_k$ defined over $k$ such that $f=f_k\otimes_kk^s$?  A
necessary condition for this to hold is that the branch locus $S$ is
$\Gamma_k$-invariant. We will assume this from now on, and let
$S_k\subset\PP^1_k$ be the closed subset of $\PP^1_k$ corresponding to
$S$ and $U_k:=\PP^1_k-S_k$ its complement. By definition, $k$ is a
field of definition of $f$ if and only if $\Phi$ can be extended to a
morphism $\Phi_k:\pi_1^t(U_k)\to G$:
\[\renewcommand{\arraystretch}{1.5}
  \begin{array}{ccccccccc}
     1 &\to& \pi_1^t(U) &\To& \pi_1^t(U_k) &\To& \Gamma_k &\to& 1 \\
          && \rupfeil{\Phi} & \;\;\swarrow{\scriptstyle\Phi_k} &&&&&  \\
          && G  &&&&&&  \\
  \end{array}
\]  
Let us fix a section $s:\Gamma_k\to\pi_1^t(U_k)$ of the natural
projection $\pi_1^t(U_k)\to\Gamma_k$, and let $\Gamma_k$ act on
$\pi_1^t(U)$ as follows:
$\op{\sigma}{\gamma}:=s(\sigma)\,\gamma\,s(\sigma)^{-1}$.

\begin{defn} \label{fomdef}
  We say that $k$ is a {\em field of moduli} of the $G$-cover 
  $f:Y\to\PP^1_{k^s}$ if for all $\sigma\in\Gamma_k$ there
  exists an element $h_\sigma\in G$ such that
  \[
       \Phi(\op{\sigma}{\gamma}\;) \;=\; 
          h_\sigma\,\Phi(\gamma)\,h_\sigma^{-1}, \qquad
          \mbox{for all\ } \gamma\in\pi_1^t(U).
  \]
\end{defn}

Note that the field of moduli condition does not depend on the choice
of the section $s$. Obviously, in order for $k$ to be a field of
definition of $f$, it is necessary that $k$ be a field of moduli.
However, this is not a sufficient condition, in general. Let $C$ be
the center of $G$ and $\Gb:=G/C$ the quotient. Assume that $k$ is a
field of moduli of $f$, let $h_\sigma$, for each $\sigma\in\Gamma_k$,
be as in Definition \ref{fomdef} and set $\phib(\sigma):=\hb_\sigma$
(the class of $h_\sigma$ in $\Gb$). Clearly, $\phib:\Gamma_k\to\Gb$ is
well defined and a group homomorphism. It follows that for all
$\sigma,\tau\in\Gamma_k$, the element
\[
     c_{\sigma,\tau} := h_\sigma\,h_\tau\,h_{\sigma\tau}^{-1}
\]
lies in the center $C$, and that $(\sigma,\tau)\mapsto
c_{\sigma,\tau}$ is a $2$-cocycle (where $C$ is regarded as constant
$\Gamma_k$-module). Let $\omega\in H^2(k,C)$ be the cohomology class
of this cocycle. By definition, $\omega$ is the obstruction for the
existence of a (weak) solution $\phi$ of the central embedding problem
$(q,\phib)$, see e.g.\ \cite[IV.6.1]{MM} (here $q:G\to\Gb$ is the
natural map).
\begin{equation} \label{embeq}
  \renewcommand{\arraystretch}{1.5}
  \begin{array}{ccccccccc}
       &&&&&& \Gamma_k       && \\ 
       &&&&& \!\!\!\!\!\!\op{\phi}{\swarrow}& \rupfeil{\phib} && \\  
     1 &\to& C &\To& G &\lpfeil{q}    & \Gb &\to& 1. \\
  \end{array}
\end{equation}
A formal verification shows (see e.g.\ \cite{DebDou1}):

\begin{prop} \label{fodprop}
  The class $\omega$ is independent of the choice of the section
  $s:\Gamma_k\to\pi_1^t(U_k)$. Moreover, $\omega=0$ if and only if $k$
  is a field of definition of $f$.
\end{prop}

We call $\omega$ the {\em f.o.d.-obstruction} for the $G$-cover $f$,
relative to $k$ (`f.o.d.' stands for `field of definition'). The rest
of the paper is concerned with studying and, if possible, computing
$\omega$, in various situations.  First we make some general remarks.

\begin{rem} \label{fodrem1}
  \begin{enumerate}
  \item
    The formation of $\omega$ is functorial in $k$. More precisely, let
    $K/k$ be a field extension,
    $f_{K^s}:Y\otimes_{k^s}K^s\to\PP^1_{K^s}$ the base change of $f$ to
    $K^s$ (we embed $k^s$ in $K^s$) and $\omega|_K$ the image of
    $\omega$ under the restriction homomorphism $H^2(k,C)\to
    H^2(K,C)$; then $\omega|_K$ is the f.o.d.-obstruction for
    $f_{K^s}$. 
  \item A $G$-cover $f:Y\to \PP^1_{k^s}$ with $r$ branch points
    corresponds to a $k^s$-point $[f]:\Spec k^s\to\Hu_r(G)$ on a
    certain moduli scheme $\Hu_r(G)$, called the {\em Hurwitz space}.
    The field $k$ is a field of moduli for $f$ if and only if $[f]$ is
    $k$-rational, i.e. factors through a morphism $\Spec k\to\Hu_r(G)$.
    See \cite{FriedVoe91} and \cite{diss}.
  \item There exists a cohomology class $\tilde{\omega}\in H^2_{\rm
      \acute{e}t}(\Hu_r(G),C)$ which specializes to $\omega$, i.e.\ 
    $\omega=[f]^*\tilde{\omega}$. The class $\tilde{\omega}$
    represents the obstruction for the existence of a global versal
    $G$-cover over $\Hu_r(G)$, see \cite{diss} or \cite{DebDouEms00}.
  \end{enumerate}
\end{rem}

For a study of \ f.o.d.-obstructions for more general (e.g.\
non-Galois) covers, see \cite{DebDou1} and \cite{DebDouEms00}.


\subsection{The Hurwitz description}
\label{hurw}

We let $k$, $S=\{z_1,\ldots,z_r\}$, $U:=\PP^1_{k^s}-S$ and $U_k$ be as
before.  We assume that $k$ has characteristic $0$; we may therefore
write $\pi_1(U)$ instead of $\pi_1^t(U)$. Let us choose a $k$-rational
{\em base point} $z_0$ on $U$. The point $z_0$ will either be a
$k$-rational point $\Spec k\to U_k$ or a ``tangential base point''
$\Spec k((z))\to U_k$, see \cite{Ihara94}. In both cases, we obtain a
section $s:\Gamma_k\to\pi_1(U_k)$ (unique up to an inner automorphism
of $\pi_1(U)$) and hence an action of $\Gamma_k$ on $\pi_1(U)$. We
will write $\pi_1(U,z_0)$ to denote the profinite group $\pi_1(U)$
together with this $\Gamma_k$-action.

We denote by $\ZZh(1):=\varprojlim_n \mmu_n$ the Tate module of
$\mathbb{G}_m$. Let 
\begin{equation} 
   \Pi \;:=\; \gen{\gamma_1,\ldots,\gamma_r\mid \prod_i\gamma_i=1}
\end{equation}   
be the free profinite group with generators
$\gamma_1,\ldots,\gamma_r$, subject to the usual product one relation.

\begin{defn} \label{repdef}
  An isomorphism
  $\rho:\Pi\iso\pi_1(U,z_0)$ of profinite groups is a {\em
presentation\/} of $\pi_1(U,z_0)$ if it has these 
  properties:
  \begin{enumerate}
  \item $\rho(\gamma_i)$ generates an inertia subgroup
    $I_{z_i}\subset\pi_1(U,z_0)$ corresponding to a point $z_i\in S$,
    and
  \item under the natural identification $I_{z_i}\cong\ZZh(1)$, all
    the $\rho(\gamma_i)$ correspond to the same element in 
    $\ZZh(1)$.
  \end{enumerate}
\end{defn}

It is a well known fact that such a presentation $\rho$ always exists;
let us choose one. To simplify the notation, we will usually identify
$\gamma_i$ with $\rho(\gamma_i)$. The {\em branch cycle argument}
states that for all $\sigma\in\Gamma_k$ there exist elements
$\beta_1,\ldots,\beta_r\in\pi_1(U,z_0)$ such that
\begin{equation} \label{hurweq1}
    \op{\sigma}{\gamma_i} \;=\; 
  \beta_i\;\gamma_{\sigma(i)}^{\chi(\sigma)}\,\beta_i^{-1}.
\end{equation}
Here $\chi:\Gamma_k\to\ZZh^\times$ is the cyclotomic character and
$\sigma(i)$ is defined by $\sigma(z_i)=z_{\sigma(i)}$.

Let $f:Y\to\PP^1_{k^s}$ be a $G$-cover with branch locus $S$,
corresponding to a homomorphism $\Phi:\pi_1(U,z_0)\to G$. Setting
$g_i:=\Phi(\gamma_i)$, we obtain an $r$-tuple $\g=(g_1,\ldots,g_r)$ of
generators of $G$ such that $\prod_ig_i=1$.  We write $\Ni_r(G)$ for
the set of all such $\g$ and $\Ni_r\inn(G):=\Ni_r(G)/\Inn(G)$ for the
quotient of this set under the action of inner automorphisms of $G$,
compare \cite{FriedVoe91}.  We obtain a bijection between
$\Ni_r\inn(G)$ (the set of {\em Nielsen classes} of length $r$) and
the set of isomorphism classes of $G$-covers with branch locus $S$. We
say that the $G$-cover $f$ has {\em Hurwitz description}
$[\g]\in\Ni_r\inn(G)$ (with respect to the presentation $\rho$).

The $\Gamma_k$-action on $\pi_1(U,z_0)$ induces a $\Gamma_k$-action on
$\Hom(\pi_1(U,z_0),G)$, and therefore on $\Ni_r(G)$. Our convention is
to let $\Gamma_k$ act on $\Ni_r(G)$ from the {\em right}, i.e. we define
\[
    \g^\sigma \;:=\; 
      (\tilde{g}_1,\ldots,\tilde{g}_r),\quad\text{with}\quad
    \tilde{g}_i := \Phi(\sigma\,\gamma_i\,\sigma^{-1}).
\]
We remark that this is not a standard convention in the literature.
The branch cycle argument \eqref{hurweq1} becomes
\[
      \tilde{g}_i \;=\; 
    b_i\;g_{\sigma(i)}^{\chi(\sigma)}\,b_i^{-1},
\]
with $b_i:=\Phi(\beta_i)$. The $\Gamma_k$-action on $\Ni_r(G)$
induces a $\Gamma_k$-action on $\Ni_r\inn(G)$. The following
Proposition merely rephrases the definitions of Section \ref{fod}.

\begin{prop} \label{hurwprop}
  Let $f:Y\to\PP^1_{k^s}$ be a $G$-cover with branch locus $S$ and
  Hurwitz description $\g$ (with respect to a presentation
  $\rho:\Pi\iso\pi_1(U,z_0)$). Then
  \begin{enumerate}
  \item 
    $k$ is a field of moduli of $f$ if and only if
    $[\g]^\sigma=[\g]$ for all $\sigma\in\Gamma_k$.
  \item
    If (i) holds then there is a unique homomorphism
    $\phib:\Gamma_k\to\Gb$ such that 
    \[
          \g^\sigma \;=\; h_\sigma\,\g\,h_\sigma^{-1},
    \]
    where $h_\sigma\in G$ is a lift of
    $\phib(\sigma)\in\Gb$, for all $\sigma\in\Gamma_k$. 
  \item 
    $k$ is a field of definition of $f$ if and only if it is a field
    of moduli and the homomorphism $\phib$ in (ii) lifts to a
    homomorphism $\phi:\Gamma_k\to G$.
  \end{enumerate}
\end{prop}

The formula $\g^\sigma = h_\sigma\,\g\,h_\sigma^{-1}$ may seem
counterintuitive at first. But we have to keep in mind that the
homomorphism $\phib$ depends on $\g$, and not only on its Nielsen
class $[\g]$. In more geometric terms, $\g$ corresponds to a {\em
  pointed $G$-cover} $f:(Y,y_0)\to(\PP^1,z_0)$. If $k$ is a field of
definition of $f$ and $\phi:\Gamma_k\to G$ is as in Proposition
\ref{hurwprop} (iii), then there exists a unique model
$f_k:Y\to\PP^1_k$ of $f$ such that the fiber $f^{-1}(z_0)$, as a
$G$-torsor over $k$, corresponds to $\phi\in H^1(k,G)$.

\begin{exa} \label{tangexa}
  Let $r:=3$, $k:=\QQ$, and let $\ell$ be an odd prime. We set
  $\ell^*:=(-1)^{(\ell-1)/2}\ell$,
  $S:=\{0,\sqrt{\ell^*},-\sqrt{\ell^*}\}\subset\PP^1_{\QQb}$ and
  $U:=\PP^1_{\QQb}-S$. We identify $\QQ(z)$ with the function field of
  $U_{\QQ}:=\PP^1_{\QQ}-S$ and let $z_0:\Spec\QQ((z))\to
  U_{\QQ}$ be the tangential base point at $0$, with parameter $z$.
  
  It is clear from \cite{Ihara94} that there exists a presentation
  $\rho:\Pi\iso\pi_1(U,z_0)$ such that $\gamma_1\in\pi_1(U,z_0)$
  corresponds to the closed path $t\mapsto e^{2\pi i t}$ (for $0\leq
  t\leq1$). Then
  \begin{equation} \label{tangeq2}
    \op{\sigma}{\gamma_1} \;=\; \gamma_1^{\chi(\sigma)}, \qquad
    \op{\sigma}{\gamma_2} \;\sim\; \gamma_{\sigma(2)}, \qquad
    \op{\sigma}{\gamma_3} \;\sim\; \gamma_{\sigma(3)}, 
  \end{equation} 
  for all $\sigma\in\Gamma_{\QQ}$. Here $\sigma$ fixes the indices $2$
  and $3$ if $\op{\sigma}{\sqrt{\ell^*}}=\sqrt{\ell^*}$ and permutes
  them if $\op{\sigma}{\sqrt{\ell^*}}=-\sqrt{\ell^*}$.
  
  From now on we assume that $\ell\not\equiv \pm 1 \pmod{8}$, and we
  set $G:={\rm SL}_2(\ell)$. Let $\ell A$, $\ell B$ be the two
  conjugacy classes of $G$ containing elements of order $\ell$ and
  $4A$ the unique conjugacy class containing elements of order $4$.
  The classes $\ell A$ and $\ell B$ are conjugate over the quadratic
  extension $\QQ(\sqrt{\ell^*})/\QQ$. By \cite{Voelklein} I.3.3.6, the
  triple $(4A,\ell A,\ell B)$ is {\em rigid}, i.e.\ there exists
  exactly one class $[g_1,g_2,g_3]\in\Ni_3\inn(G)$ such that $g_1\in
  4A$, $g_2\in\ell A$ and $g_3\in\ell B$. Let $f:Y\to\PP^1_{\QQ}$ be
  the $G$-cover with Hurwitz description $[\g]=[g_1,g_2,g_2]$ (with
  respect to $\rho$). By Rigidity and \eqref{tangeq2},
  $[\g]^\sigma=[\g]$ for all $\sigma\in\Gamma_{\QQ}$. Therefore,
  $\QQ$ is a field of moduli of $f$.
  
  The center $C$ of $G$ consists of the two diagonal matrices $\pm I$;
  we identify it with $\{\pm 1\}$. Let $\omega\in H^2(\QQ,\pm 1)$ be
  the f.o.d.-obstruction of $f$. We may identify $H^2(\QQ,\pm 1)$ with
  $\Br_2(\QQ)$, the $2$-torsion of the Brauer group of $\QQ$.

\begin{prop} \label{2llprop}
  We have $\omega=(-1,-1)$, i.e.\ $\omega$ is represented (as an
  element of $\Br_2(\QQ)$) by the quaternion algebra $\QQ[i,j\mid
  i^2=-1,\,j^2=-1,\,ij=-ji]$. In particular, $\QQ$ is not a field of
  definition of $f$.
\end{prop}

It follows that the regular ${\rm PSL}_2(\ell)$-extensions of $\QQ(z)$
with ramification of type $(2,\ell,\ell)$ (see e.g.\
\cite{SerreTopics}, Section 8.3.3) do not lift to an ${\rm
  SL}_2(\ell)$-extension. We will see in Section \ref{real} that in
fact no ${\rm PSL}_2(\ell)$-extension with three branch points lifts to
an ${\rm SL}_2(\ell)$-extension.

\begin{proof}
  The class $\omega\in H^2(\QQ,\pm 1)$ is the obstruction for lifting
  the homomorphism $\phib:\Gamma_{\QQ}\to\Gb={\rm PSL}_2(\ell)$ to a
  homomorphism $\phi:\Gamma_{\QQ}\to G={\rm SL}_2(\ell)$. By the
  definition of $\phib$ and by Equation \eqref{tangeq2}, we have (by
  abuse of notation) \begin{equation} \label{tangeq3}
  \phib(\sigma)\,g_1\,\phib(\sigma)^{-1} \;=\; g_1^{\,\chi(\sigma)},
  \qquad \text{for all}\;\; \sigma\in\Gamma_{\QQ}.  \end{equation} In
  particular, the image of $\phib$ is contained in $\bar{N}:=N/C$,
  where $N$ is the normalizer in $G$ of the cyclic subgroup generated
  by $g_1$. It is shown e.g.\ in \cite{Huppert}, Abschnitt II.8, that
  $\bar{N}$ is a dihedral group of order $4n$, where $n:=(\ell-1)/4$
  if $\ell\equiv 1\pmod{4}$ and $n:=(\ell+1)/4$ otherwise. Let
  $\bar{H}\lhd\bar{N}$ be the unique cyclic normal subgroup of order
  $n$. By our assumption on $\ell$, $n$ is odd. Therefore, there
  exists a unique cyclic normal subgroup $H\lhd N$ of order $n$
  mapping onto $\bar{H}$. Since $\bar{N}/\bar{H}\cong\ZZ/2\times\ZZ/2$
  and every involution in $\Gb$ lifts to an element of order $4$ in
  $G$, $N/H$ is isomorphic to the quaternion group $Q_8$. It is easy
  to check that in the following diagram 
  \begin{equation} 
  \label{tangeq4} \renewcommand{\arraystretch}{1.5}
  \begin{array}{ccccc} 
    \{\pm 1\} & \To & N & \To & \bar{N} \\
    \big\downarrow & & \big\downarrow & & \big\downarrow \\ 
    \{\pm 1\} &\To & N/H\cong Q_8 & \To & N/CH\cong\ZZ/2\times\ZZ/2 \\
  \end{array}
  \end{equation}   
  the group extension in the upper row is the pullback of the lower
  row with respect to the projection $\bar{N}\to N/CH$. Let
  $\psib:\Gamma_{\QQ}\to N/CH$ be the composition of $\phib$ with the
  projection $\bar{N}\to N/CH$. It follows from \cite{Brown}, Exercise
  IV.3.1, that $\omega=\psib^*\eta$, where $\eta\in H^2(N/CH,C)$ is
  the cohomology class corresponding to the lower row of
  \eqref{tangeq4}. In other words, $\omega$ is the
  obstruction for lifting $\psib$ to a homomorphism
  $\psi:\Gamma_{\QQ}\to N/H\cong Q_8$. Therefore, by a result of Witt
  (see \cite{Froehlich85}, (7.7) (ii)), 
  \[ 
     \omega = (-1,a) + (-1,b) + (a,b) \in \Br_2(\QQ), 
  \] 
  where $E=\QQ[x,y\mid x^2=a,\,y^2=b]$, $a,b\in\QQ$, is
  the \'etale $\QQ$-algebra corresponding to $\psib$. Let
  $\bar{H}'\lhd\bar{N}$ be the maximal cyclic normal subgroup, of
  order $2n$. By \eqref{tangeq3}, the composition of $\psib$ with
  $N/CH\to\bar{N}/\bar{H}'\cong\ZZ/2\cong(\ZZ/4)^\times$ equals the
  cyclotomic character modulo $4$. Hence we may choose $a=-1$, and
  $\omega=(-1,-1)$ follows.
\end{proof}
\end{exa}


\subsection{$G$-covers over the real numbers} \label{real}

In this section we discuss the case $k=\RR$. Since the action of
$\Gamma_{\RR}$ on $\pi_1(U,z_0)$ is known, we get some easy results
and examples which illustrate several key concepts of this paper. For
more general results on related questions, see
\cite{BaileyFriedTV}, in particular \S 6.

Let $f:Y\to\PP^1_{\CC}$ be a $G$-cover defined over the complex
numbers, and assume that the branch locus $S$ of $f$ is defined over
$\RR$. For simplicity, we will also assume that $\infty\not\in S$. We
can write $S=\{z_1,\ldots,z_r\}$, with $\bar{z}_i=z_{2s-i+1}$ for
$i=1,\ldots,s$, and $z_{2s+j}\in\RR$ for $j=1,\ldots,r-2s$. For
$z,w\in\CC$, we write $z\succ w$ if either $\Im z>\Im w$ or $\Im z=\Im
w$ and $\Re z >\Re w$. In this notation, we may assume that
$z_1\succ\ldots\succ z_s$ and $z_{2s+1}<\ldots<z_r$. Let
$U:=\PP^1_{\CC}-S$ and choose a real base point $z_0<z_{2s+1}$. By
\cite[Section I.1.1]{MM}, there exists a presentation
$\rho:\Pi\iso\pi_1(U,z_0)$ with the following properties: (i)
$\gamma_i$ is represented by a simple loop based at $z_0$ winding
counterclockwise around $z_i$, (ii) these loops are pairwise disjoint
and (iii) the action of complex conjugation $\kappa\in\Gamma_{\RR}$ is
given by 
\begin{alignat}{2} \label{realeq2}
   \op{\kappa}{\gamma_i} &\;=\; \gamma_{2s-i+1}^{-1}, &\qquad
      \text{for}\;\; & i=1,\ldots,2s,\\
\intertext{and} \label{realeq3}
   \op{\kappa}{\gamma_{2s+j}} &\;=\; 
         \beta_j\,\gamma_{2s+j}^{-1}\,\beta_j^{-1}, &\qquad
        \text{for}\;\;  &j=1,\ldots,r-2s,
\end{alignat}
where $\beta_j:=\gamma_{2s+1}\cdot\ldots\cdot\gamma_{2s+j-1}$ (see
\cite[Fig.\ 1.2]{MM}). We say that $\rho$ is a {\em normalized
  presentation} of $\pi_1(U,z_0)$.

\begin{prop} \label{realprop}
  Let $\g=(g_1,\ldots,g_r)$ be the Hurwitz description for the
  $G$-cover $f:Y\to\PP^1_{\CC}$, with respect to a normalized
  presentation $\rho$, as above. Then
  \begin{enumerate}
  \item
    $\RR$ is a field of moduli for $f$ if and only if there exists
    $b\in G$ such that
    \[
       b\;g_i\,b^{-1} \;=\; g_{2s-i+1}^{-1}
    \]
    for $i=1,\ldots,2s$, and 
    \[
       b\;g_{2s+j}\,b^{-1} \;=\; h_j\,g_{2s+j}^{-1}\,h_j^{-1}
    \]
    for $j=1,\ldots,r-2s$, with $h_j:=g_{2s+1}\cdot\ldots\cdot
    g_{2s+j-1}$.
  \item
    $\RR$ is a field of definition of $f$ if and only there exists
    $b$ as in {\rm (i)} such that $b^2=1$.
  \end{enumerate}
\end{prop}

\begin{proof}
  Follows from Proposition \ref{hurwprop} and equations
  \eqref{realeq2} and \eqref{realeq3}.
\end{proof}

For certain groups $G$, the condition in Proposition \ref{realprop}
(ii) is rather restrictive, as shown by the next theorem.

\begin{thm} \label{realthm}
  Let $G$ be a finite group. Assume that all elements of $G$ of order
  $2$ lie in the center of $G$. Let $f:Y\to\PP^1_{\CC}$ be a
  $G$-cover, with branch locus $S$ defined over $\RR$, and let $\g$ be
  a Hurwitz description of $f$ with respect to a normalized
  presentation $\rho:\Pi\to\pi_1(U,z_0)$ (see above). Then $\RR$ is a
  field of definition of $f$ if and only if
  \begin{equation} \label{realeq4}
     \g = (\,g_1,\ldots,g_s,g_s^{-1},\ldots,g_1^{-1},
            \,c_1,\ldots,c_{r-2s}),
  \end{equation}
  with $c_j^2=1$ and $\prod_j c_j=1$.
\end{thm}

The
quaternion group $Q_8$ and ${\rm SL}_2(\ell^n)$, for $\ell\not=2$ verify
the conditions of the theorem. Elements of order 2 in $A_n$ are products
of $2s$ disjoint 2-cycles. They lift to have order 4 in the $2$-fold 
coverings $\tilde{A}_n$ of $A_n$ exactly when $s$ is odd
\cite[Prop.~5.8]{BaileyFriedTV}. So this also holds for $\tilde{A}_n$ if
$4\le n\le 7$.

\begin{proof}
  If $\g$ is as in \eqref{realeq4}, then $\op{\kappa}{\g}=\g$, so $f$
  can be defined over $\RR$. Conversely, if $f$ can be defined over
  $\RR$ then there exists $b\in G$ of order $2$ such that
  $b\,g_i\,b^{-1}=g_{2s-i+1}^{-1}$ and
  $b\,g_{2s+j}\,b^{-1}=h_j\,g_{2s+j}^{-1}\,h_j^{-1}$, by
  Proposition \ref{realprop} (ii). By the hypothesis of the theorem,
  $b$ is a central element of $G$. Therefore,
  $g_i=g_{2s-i+1}^{-1}$, for $i=1,\ldots,s$. Moreover,
  $g_{2s+1}=g_{2s+1}^{-1}$, so $c_1:=g_{2s+1}$ has order $2$ and hence
  lies in the center of
  $G$. Continuing by induction, we find that $c_j:=g_{2s+j}$ is
  central of order $2$, for $j=1,\ldots,r-2s$. This proves the
  theorem.
\end{proof}

\begin{exa} \label{realexa}
  Let $G:=\At_5$ be the nonsplit central extension of $A_5$ by $\pm
  1$. Let $g_1,g_2\in G$ be the unique lifts of order $3$ of
  $\gb_1:=(1\,2\,3),\;\gb_2:=(3\,4\,5)\in\Gb=A_5$. Set
  $\g:=(\,g_1,\,g_2,\,g_2^{-1},g_1^{-1})\in\Ni_4(G)$. For $t\in\CC$,
  set
\[
    S_t \;:=\; \{\;z\in\CC \mid z^4+(2t-2)z^2+t^2+2t+1=0\;\}
        \;=\;  \{\;\pm1 \pm\sqrt{-t}\;\}
\]
and $U_t:=\PP^1_{\CC}-S_t$. Let us fix a value $t_0\in\RR$ such that
$0<t_0<1$, and let $\rho_{t_0}$ be the normalized presentation of
$\pi_1(U_{t_0},0)$, as shown on the left hand side of Figure
\ref{realfig}. We let $f_{t_0}:Y_{t_0}\to\PP^1_{\CC}$ be the $G$-cover
with Hurwitz description $\g$, with respect to $\rho_{t_0}$. By
Theorem \ref{realthm}, $f_{t_0}$ can be defined over $\RR$.

\begin{figure}[bt]
\begin{minipage}{6cm}
  \begin{center}
    {\setlength{\unitlength}{0.00045in}

\begin{picture}(4962,2289)(0,-10)
\put(1350,2037){\circle*{90}}
\put(4050,2037){\circle*{90}}
\put(1350,237){\circle*{90}}
\put(4050,237){\circle*{90}}
\put(1350,2037){\circle{450}}
\put(4050,2037){\circle{450}}
\put(1350,237){\circle{450}}
\put(4050,237){\circle{458}}
\drawline(450,1137)(4950,1137)
\drawline(2700,1137)(1530,1902)(1530,1902)
\drawline(2700,1137)(3870,1902)
\drawline(2700,1137)(1530,372)
\drawline(2700,1137)(3870,372)
\drawline(1395,2262)(1305,2262)
\drawline(1425.000,2292.000)(1305.000,2262.000)(1425.000,2232.000)
\drawline(4095,2262)(4005,2262)
\drawline(4125.000,2292.000)(4005.000,2262.000)(4125.000,2232.000)
\drawline(1305,12)(1395,12)
\drawline(1275.000,-18.000)(1395.000,12.000)(1275.000,42.000)(1275.000,-18.000)
\drawline(4005,12)(4095,12)
\drawline(3975.000,-18.000)(4095.000,12.000)(3975.000,42.000)(3975.000,-18.000)
\put(0,1050){$\RR$}
\put(4500,1947){$\gamma_1$}
\put(700,2000){$\gamma_2$}
\put(700,100){$\gamma_3$}
\put(4500,147){$\gamma_4$}
\end{picture}
}

  \end{center}
\end{minipage}
\hfill
\begin{minipage}{6cm}
   \begin{center}
     {\setlength{\unitlength}{0.00045in}
\begin{picture}(4962,1275)(0,70)
\drawline(3420.000,891.000)(3351.722,914.132)(3281.251,929.324)
   (3209.508,936.377)(3137.427,935.201)(3065.952,925.809)
   (2996.015,908.324)(2928.529,882.976)(2864.375,850.094)
   (2804.390,810.109)(2749.359,763.541)(2700.000,711.000)
\drawline(1980.000,576.000)(2033.830,527.335)(2095.015,488.316)
   (2161.844,460.033)(2232.450,443.276)(2304.861,438.514)
   (2377.054,445.880)(2447.011,465.167)(2512.778,495.838)
   (2572.518,537.034)(2624.562,587.606)(2667.455,646.140)
   (2699.999,711.000)
\drawline(2700.000,711.000)(2745.748,655.468)(2794.538,602.588)
   (2846.217,552.527)(2900.621,505.442)(2957.579,461.482)
   (3016.913,420.784)(3078.435,383.478)(3141.952,349.680)
   (3207.265,319.496)(3274.167,293.022)(3342.448,270.341)
   (3411.894,251.524)(3482.285,236.631)(3553.401,225.708)
   (3625.017,218.789)(3696.909,215.897)(3768.849,217.041)
   (3840.613,222.216)(3911.973,231.408)(3982.706,244.586)
   (4052.588,261.709)(4121.401,282.723)(4188.926,307.563)
   (4254.953,336.150)(4319.273,368.394)(4381.684,404.194)
   (4441.990,443.437)(4500.000,486.000)
\drawline(2700.000,711.000)(2654.252,766.532)(2605.462,819.412)
   (2553.783,869.473)(2499.379,916.558)(2442.421,960.518)
   (2383.087,1001.216)(2321.565,1038.522)(2258.048,1072.320)
   (2192.735,1102.504)(2125.833,1128.978)(2057.552,1151.659)
   (1988.106,1170.476)(1917.715,1185.369)(1846.599,1196.292)
   (1774.983,1203.211)(1703.091,1206.103)(1631.151,1204.959)
   (1559.387,1199.784)(1488.027,1190.592)(1417.294,1177.414)
   (1347.412,1160.291)(1278.599,1139.277)(1211.074,1114.437)
   (1145.047,1085.850)(1080.727,1053.606)(1018.316,1017.806)
   (958.010,978.563)(900.000,936.000)
\put(900,711){\circle*{90}}
\put(4500,711){\circle*{90}}
\put(1800,711){\circle*{90}}
\put(3600,711){\circle*{90}}
\put(900,711){\circle{450}}
\put(1800,711){\circle{450}}
\put(3600,711){\circle{450}}
\put(4500,711){\circle{450}}
\drawline(450,711)(4950,711)
\drawline(3555,486)(3645,486)
\drawline(3525.000,456.000)(3645.000,486.000)(3525.000,516.000)(3525.000,456.000)
\drawline(4545,936)(4455,936)
\drawline(4575.000,966.000)(4455.000,936.000)(4575.000,906.000)(4575.000,966.000)
\drawline(1845,936)(1755,936)
\drawline(1875.000,966.000)(1755.000,936.000)(1875.000,906.000)(1875.000,966.000)
\drawline(855,486)(945,486)
\drawline(825.000,456.000)(945.000,486.000)(825.000,516.000)(825.000,456.000)
\put(0,576){$\RR$}
\put(3465,1250){$\gamma_1$}
\put(4410,1250){$\gamma_4$}
\put(810,81){$\gamma_2$}
\put(1710,36){$\gamma_3$}
\end{picture}
}

   \end{center}
\end{minipage}
\caption{}
\label{realfig}
\end{figure}

Let $\theta:[0,1]\to \CC$ be the path $s\mapsto\theta_s:=t_0\,e^{\pi
  is}$, and let $f_{-t_0}:Y_{-t_0}\to\PP^1_{\CC}$ be the $\At_5$-cover
with branch locus $S_{-t_0}$ obtained from $f_{t_0}$ by ``analytic
continuation of the branch locus along $s\mapsto S_{\theta_s}$'' (see
e.g.\ \cite[Section II.10.1]{Voelklein} for precise definitions).
Thus, $f_{-t_0}$ has branch locus $S_{-t_0}=\{\pm 1\pm\sqrt{t_0}\}$,
and Hurwitz description $\g$ with respect to the presentation
$\rho_{-t_0}$, as shown on the right hand side of Figure
\ref{realfig}. Note that the presentation $\rho_{-t_0}$ is not
normalized. From Figure \ref{realfig}, we see that the action of
complex conjugation on $\pi_1(U_{-t_0},0)$, with respect to
$\rho_{-t_0}$, is given by
\begin{equation} \label{realeq5}
   \op{\kappa}{\gamma_1}=\gamma_1^{-1}, \quad
   \op{\kappa}{\gamma_2}=\gamma_3^{-1}\,\gamma_2^{-1}\,\gamma_3,\quad
   \op{\kappa}{\gamma_3}=\gamma_3^{-1}, \quad
   \op{\kappa}{\gamma_4}=\gamma_1^{-1}\,\gamma_4^{-1}\,\gamma_1.
\end{equation}
A short computation reveals how this new action affects the tuple
$\g$:
\begin{equation} \label{realeq6}
   \g^\kappa \;=\; (\,g_1^{-1},g_2^{-1},\,g_2,\,g_1\,) 
                   \;=\; b\,\g\,b^{-1},
\end{equation}
where $b\in\At_5$ is a lift of $(1\,2)(4\,5)\in \Gb=A_5$. Therefore,
$\RR$ is a field of moduli of $f_{-t_0}$. On the other hand, any
element $b\in\At_5$ such that \eqref{realeq6} holds has order $4$
(there are exactly two of them), so $\RR$ is not a field of definition
of $f_{-t_0}$. In other words: for $t\in(-1,1)-\{0\}$, let
$\omega_t\in H^2(\RR,\{\pm 1\})$ be the f.o.d.-obstruction for
$f_t$. Then $\omega=0$ for $t>0$ and $\omega=(-1,-1)\not=0$ for $t<0$. 
\end{exa}

Example \ref{realexa} illustrates a key point of this paper: the
obstruction $\omega_t$, as a function of $t\in (-1,1)\subset\RR$, has
a `jump' at $t=0$. In the next section we will give the following
algebraic interpretation for this phenomenon.  Let
$\tilde{f}:\tilde{Y}\to\PP^1_{k^s}$ be the $G$-cover over $k^s$, with
$k=\RR((t))$, corresponding to the family $f_t$ in an infinitesimal
neighborhood of $t=0$. Then $k$ is a field of moduli for $\tilde{f}$,
and the corresponding f.o.d.-obstruction $\omega\in H^2(k,\{\pm 1\})$
has a {\em pole} at $t=0$, i.e.\ the {\em residue} $r_v(\omega)\in
H^1(\RR,\{\pm 1\})$ of $\omega$ at the valuation $v$ of $k$ does not
vanish.



\section{Galois covers over henselian fields}

In Section \ref{residue} we recall the definition of the residue map
and we state a theorem that says that the residue of the
f.o.d.-obstruction of a $G$-cover vanishes if the cover has good
reduction. In Section \ref{comp} we define {\em compatible
automorphisms} of the free profinite group $\Pi$, with respect to an
{\em ordered tree}. This is preparatory work for Section \ref{adm},
where we state some results on the Galois action on $\pi_1(U)$. Here
we restrict our attention to henselian fields with residue
characteristic $0$. In Section \ref{compute} we explicitly compute the
f.o.d.-obstruction for some $\At_5$-covers over finite extensions of
$\QQ((t))$. In Section \ref{spec} we explain how to `specialize'
our previous results to $p$-adic fields.

\subsection{The residue map} 
\label{residue}

Throughout this section, we assume that $k$ has characteristic $0$
and is henselian with respect to a discrete valuation $v$. We denote
the residue field of $v$ by $k_0$ and assume that $k_0$ is perfect. We
write $\Gamma_k^t$ for the maximal tame quotient of $\Gamma_k$ and
$I^t\lhd\Gamma_k^t$ for its tame inertia subgroup. We choose a section
for the natural map $\Gamma_k^t\to\Gamma_{k_0}$ and identify
$\Gamma_k^t$ with $I^t\rtimes\Gamma_{k_0}$. We remark that, as a
$\Gamma_{k_0}$-module, $I^t$ is canonically isomorphic to
$\ZZh'(1):=\varprojlim_n \mu_n(k_0^s)$, where $n$ runs over the
integers prime to the characteristic of $k_0$.

We are given a $G$-cover $f:Y\to\PP^1_{k^s}$ with branch locus $S$ and
field of moduli $k$. Let $\omega\in H^2(k,C)$ be the
f.o.d.-obstruction for $f$. We choose a $k$-rational branch point
$z_0$ on $U:=\PP^1_{k^s}-S$. As explained in Sections \ref{fod} and
\ref{hurw}, we obtain a certain homomorphism $\phib:\Gamma_k\to\Gb$,
and $\omega$ is the cohomological obstruction for lifting $\phib$ to a
homomorphism $\phi:\Gamma_k\to G$. We assume that the following
tameness condition holds.

\begin{cond}\ \label{tamecond}
\begin{enumerate}
\item
  The order of $C$ is prime to the characteristic of $k_0$.
\item
  The homomorphism $\phib:\Gamma_k\to\Gb$ is at most tamely ramified
  at $v$.
\end{enumerate}
\end{cond}

We may therefore regard $\phib$ as a homomorphism $\Gamma_k^t\to\Gb$
and $\omega$ as the obstruction for lifting $\phib$ to a homomorphism
$\phi:\Gamma_k^t\to G$ (Condition \ref{tamecond} (i) implies that
there is a natural isomorphism $H^2(k,C)\cong H^2(\Gamma_k^t,C)$). 

According to \cite{SerreCG}, Section II.A.2, and
Condition \ref{tamecond} (i), we have a short exact sequence
\begin{equation} \label{reseq}
  0 \to H^2(k_0,C) \To H^2(k,C) \lpfeil{r_v} H^1(k_0,C(-1)) \to 0
\end{equation}
(with $C(-1):=\Hom(I^t,C)$). This sequence can be deduced from the
Hochschild--Serre spectral sequence for the group extension determined
by the inclusion $I\inj\Gamma_k$. The map $r_v$ is called the {\em
  residue map}, and $r_v(\omega)$ is called the residue of $\omega$ at
$v$.  If $r_v(\omega)=0$ we say that $\omega$ is {\em regular} at $v$;
otherwise, $\omega$ is said to have a {\em pole} at $v$. If $\omega$
is regular at $v$, we can identify $\omega$ with a class in
$H^2(k_0,C)$; we denote this class by $\omega_0$ and refer to it as
the {\em value} of $\omega$ at $v$.

Proposition \ref{residueprop} below gives a description of
$r_v(\omega)$ in terms of $\phib$. To state it, we need some more
notation. Let us choose a topological generator $q_0$ of the inertia
group $I^t$. To simplify the notation, we will identify
$C(-1)=\Hom(I^t,C)$ with $C$ (as abelian groups), via $\eta\mapsto
\eta(q_0)$. We let $\Gamma_{k_0}$ act on $G$ as follows:
\begin{equation} \label{actioneq}
   \op{\sigma}{g} \;:=\; 
     \phib(\sigma)\,g\,\phib(\sigma)^{-1}, \qquad
     \sigma\in\Gamma_{k_0}\subset\Gamma_k^t.
\end{equation}
Finally, let us choose an element $a_0\in G$ lifting $\phib(q_0)\in\Gb$.
For all $\sigma\in\Gamma_{k_0}$, we find that
\begin{equation} \label{cocycleeq}
     c_\sigma:= a_0\,(\op{\sigma}{a_0})^{-\chi(\sigma)^{-1}}=
a_0\,((\op{\sigma} {a_0})^{\chi(\sigma)^{-1}})^{-1}
\end{equation}
is an element of $C$. 

\begin{prop} \label{residueprop}
  \begin{enumerate}
  \item 
    The map $\sigma\mapsto c_\sigma$ is a $1$-cocycle with
    values in $C\cong C(-1)$. Its cohomology class equals the residue
    $r_v(\omega)$.
  \item
    If $r_v(\omega)=0$, then $\omega_0\in H^2(k_0,C)$ is the
    obstruction for lifting the homomorphism
    $\phib_0:=\phib|_{\Gamma_{k_0}}$ to $G$.  
  \end{enumerate}
\end{prop}

\begin{proof}
  For each $\sigma\in\Gamma_{k_0}$, let us choose a lift $b_\sigma\in
  G$ lifting $\phib(\sigma)$. We may assume that $\sigma\mapsto
  b_\sigma$ is continuous.  
  Every element $\tau\in\Gamma_k^t$ can be written in a unique way as
  $\tau=\sigma\,q_0^i$, with $\sigma\in\Gamma_{k_0}$ and
  $i\in(\ZZh')^\times$. Therefore,
\begin{equation} \label{residueeq1}
    \phi'(\sigma\,q_0^i) := b_\sigma\,a_0^i
\end{equation}
defines a continuous, set-theoretic lift $\phi':\Gamma_k^t\to G$ of
$\phib$. By definition, the class $\omega\in H^2(k,C)$ is represented
by the $2$-cocycle
\begin{equation} \label{residueeq2}
   \omega_{\tau_1,\tau_2} := 
     \phi'(\tau_1)\,\phi'(\tau_2)\,\phi'(\tau_1\tau_2)^{-1}, \qquad
   \tau_1,\tau_2\in \Gamma_k^t.
\end{equation}
Using \eqref{residueeq1}, one checks that $\omega_{\tau_1,\tau_2}$ is
normalized (i.e.\ equals $1$ when one of the $\tau_j$ equals $1$) and
depends only on the class of $\tau_2$ in
$\Gamma_{k_0}=\Gamma_k^t/I^t$. By \cite{SerreCG}, Section II.A.1, the
residue $r_v(\omega)\in H^1(k_0,\Hom(I^t,C))$ is represented by the
$1$-cocycle
\begin{equation} \label{residueeq3}
   \sigma \;\longmapsto\; 
   (q\mapsto \omega_{q,\sigma}), \qquad
   \sigma\in\Gamma_{k_0},\;q\in I^t.
\end{equation}
Using \eqref{residueeq1}, \eqref{residueeq2} and the identification
$C\cong C(-1)=\Hom(I^t,C)$, \eqref{residueeq3} becomes
\[\begin{split}
   \sigma\;\longmapsto\; \omega_{q_0,\sigma}
 &=\phi'(q_0)\,\phi'(\sigma)\,\phi'(q_0\,\sigma)^{-1} \\
 &=\phi'(q_0)\,\phi'(\sigma)\,
        \phi'(\sigma\,q_0^{\chi(\sigma)^{-1}})^{-1}\\
 &=a_0\,(\op{\sigma}{a_0})^{-\chi(\sigma)^{-1}} \;=\; c_\sigma. \\
\end{split}\]
This proves the first part of the proposition. The second part is
obvious from the definition of the sequence \eqref{reseq}.
\end{proof}

\begin{cor} \label{residuecor}
  Assume that the order of $\phib(q_0)\in\Gb$ is prime to the order of
  $C$. Then $r_v(\omega)=0$.
\end{cor}

\begin{proof}
  It follows from our assumption that there exists a lift
  $a_0\in G$ of $\phib(q_0)$ such that $\ord a_0=\ord\phib(q_0)$ is
  prime to the order of $C$. By the definition of $c_\sigma$ we
  have
  \[
        \op{\sigma}{a_0}= (c_\sigma^{-1}\,a_0)^{\chi(\sigma)},
  \]
  for $\sigma\in\Gamma_{k_0}$. Choose an integer $n$ such that
  $n\equiv 1\pmod{\ord a_0}$ and $n\equiv 0\pmod{\ord C}$. Then 
  \[
      \op{\sigma}{a_0} = \op{\sigma}{(a_0^n)} = 
      (c_\sigma^{-n}\,a_0^n)^{\chi(\sigma)} = a_0^{\chi(\sigma)}.
  \]
  It follows that $c_\sigma=1$, for all $\sigma\in\Gamma_{k_0}$.
\end{proof}

Let us denote by $R$ the ring of integers of $k^s$. Let $S_R\subset
\PP^1_R$ be the closure of $S$ in $\PP^1_R$. We say that the $G$-cover
$f:Y\to \PP^1_{k^s}$ has {\em good reduction} if (i) $S_R$ is \'etale
over $\Spec R$ and (ii) $f$ extends to a finite morphism
$f_R:Y_R\to\PP^1_R$ which is tamely ramified along $S_R$ and \'etale
everywhere else. If $f$ has good reduction, then
$\bar{f}:\bar{Y}:=Y_R\otimes k_0^s\to\PP^1_{k_0^s}$ is a $G$-cover,
with branch locus $\bar{S}:=S_R\otimes k_0^s$.

\begin{thm} \label{redthm}
  Assume that $f$ has good reduction at $v$. Then $r_v(\omega)=0$.
  Moreover, the value $\omega_0\in H^2(k_0,C)$ of $\omega$ at $v$
  is the f.o.d.-obstruction for the $G$-cover $\bar{f}$ (for which
  $k_0$ is a field of moduli).
\end{thm}

This theorem can be seen as a consequence of Remark \ref{fodrem1}
(iii). Let us give a proof which does not rely on the results
of \cite{diss}.

\begin{proof}
  Let $\bar{U}:=\PP^1_{k_0^s}-\bar{S}$ and
  $\lambda:\pi_1(U,z_0)\to\pi^t_1(\bar{U},\bar{z}_0)$ the
  specialization morphism, defined in \cite{SGA1}. Since $f$ is
  assumed to have good reduction, $\Phi:\pi_1(U,z_0)\to G$ factors
  through $\lambda$ and the resulting morphism
  $\bar{\Phi}:\pi^t_1(\bar{U},\bar{z}_0)\to G$ corresponds to the
  reduction $\bar{f}$. Because $\lambda$ is $\Gamma_k$-equivariant and
  the inertia group $I\lhd\Gamma_k$ acts trivially on
  $\pi_1^t(\bar{U},\bar{z})$, the homomorphism $\phib:\Gamma_k\to\Gb$
  is unramified, i.e.\ corresponds to a homomorphism
  $\phib_{k_0}:\Gamma_{k_0}\to\Gb$. Thus, Proposition \ref{residueprop}
  implies that $r_v(\omega)=0$ and that $\omega_0$ is the obstruction
  for lifting $\phib_{k_0}$ to $G$, proving the theorem.
\end{proof}

Theorem \ref{redthm} above implies Theorem 3.1 of \cite{DebHar98}.
Namely, assume that the characteristic of $k_0$ does not divide the
order of $G$ and that $S_R$ is \'etale over $\Spec R$. Then $f$ has
good reduction, by \cite{Beckmann89}, and hence $r_v(\omega)=0$.  If
we assume in addition that ${\rm cd}\,k_0\leq 1$ (e.g.\ if $k_0$ is
finite) then $\omega=0$, because the sequence \eqref{reseq} is exact.


\subsection{Ordered trees and compatible automorphisms}
\label{comp}

Let $T=(V,E)$ be a {\em tree}, i.e.\ a finite and simply connected
graph. We choose a distinguished vertex $v_0\in V(T)$ and call it the
{\em root} of $T$. The choice of $v_0$ induces a natural orientation on
$T$; we represent an edge $e\in E(T)$ as an ordered pair $e=(v_1,v_2)$
of
vertices such that $v_1$ is closer to $v_0$ than $v_2$.  For $v\in
V(T)$ let $p_v=(v_0,\ldots,v)$ be the shortest path leading from $v_0$
to $v$. We define $T_v$ as the subtree of $T$ which contains all the
vertices $v'$ such that $p_{v'}$ passes through $v$.  We set
\[
     A_v \;:=\; \{\; v' \mid (v,v')\in E(T) \;\}.
\]
For each vertex $v'$ there exists a unique vertex $v=\pre(v')$ such
that $v'\in A_v$, called the {\em predecessor} of $v'$.  A vertex
$v\in V(T)-\{v_0\}$ will be called a {\em leaf} if $A_v=\emptyset$. We
write $B(T)$ for the set of leaves of $T$.

\begin{defn}
  An {\em order} on $T$ is a bijection $\psi:\{1,\ldots,r\}\iso
  B(T)$ such that for each $v\in V(T)$ the set
  \[
       I_v \;:=\; \{\;i\,\mid\; 
                     \psi(i)\in V(T_v)\;\} \subset \{1,\ldots,r\}
  \]
  is an interval, i.e.\ $I_v=\{r',\ldots,r''\}$. The triple
  $(T,v_0,\psi)$ is called an {\em ordered tree}.
\end{defn}

An order $\psi$ on $T$ induces a strict ordering of the sets $A_v$:
for $v',v''\in A_v$ we write $v'< v''$ if $i\in I_{v'}$
and $j\in I_{v''}$ implies $i<j$. Conversely, if we choose a strict
ordering on each set $A_v$ then there exists a unique order $\psi$ on
$T$ inducing these orderings. In particular, an order $\psi$ always
exists. An automorphism of the ordered tree $(T,v_0,\psi)$ is an
automorphism $\kappa:T\iso T$ such that $\kappa(v_0)=v_0$. Via $\psi$,
$\kappa$ induces a permutation of $\{1,\ldots,r\}$. Note that $\kappa$
is uniquely determined by this permutation.

Let $\Pi=\gen{\gamma_1,\ldots,\gamma_r\mid\prod\gamma_i=1}$ be the
free profinite group defined earlier, and let $(T,v_0,\psi)$ be an
ordered tree. For each $v\in V(T)$ we define an element
\[
     \gamma_v := \prod_{i\in I_v} \gamma_i
\]
and a subgroup 
\[
     \Pi_v := \gen{\gamma_{v'} \mid v'\in A_v} \subset \Pi
\]
of $\Pi$. Note that $\gamma_{v_0}=1$, $\gamma_{\psi(i)}=\gamma_i$
and that for all $v\in V(T)$ we have
\begin{equation} 
  \gamma_v = \prod_{v'\in A_v} \gamma_{v'},
\end{equation}
where the product is taken with respect to the ordering on $A_v$
induced by $\psi$. 

\begin{defn} \label{compdef}
  An automorphism $\tau:\Pi\iso\Pi$ is {\em compatible} with
  $(T,v_0,\psi)$ if there exists an automorphism $\kappa$ of
  $(T,v_0,\psi)$, an element $\chi\in\ZZh^\times$ and elements
  $\alpha_v\in \Pi_{\pre(v)}$, for all $v\in V(T)$, such that the
  following holds. For each $v\in V(T)-\{v_0\}$, let
  $\beta_v:=\alpha_{v_1}\cdot\ldots\cdot\alpha_v$, where
  $p_v=(v_0,v_1,\ldots,v)$ is the shortest path from $v_0$ to $v$.
  Then
  \[
        \tau(\gamma_v) = 
     \beta_{\kappa(v)}\, \gamma_{\kappa(v)}^\chi\,
\beta_{\kappa(v)}^{-1}.
  \]
\end{defn}

It is easy to see that $\chi=\chi(\tau)$ and $\kappa=\kappa(\tau)$ are
uniquely determined by $\tau$.

\begin{exa} \label{compexa1}
  Let $\rho:\Pi\iso\pi_1(U,z_0)$ be a presentation, as in Definition
  \ref{repdef}. An element $\sigma\in\Gamma_k$ induces an automorphism
  of $\pi_1(U)$, and hence, via $\rho$, an automorphism of $\Pi$.
  Let $V:=\{0,1,\ldots,r\}$, $E:=\{(0,1),\ldots,(0,r)\}$ and
  $T:=(V,E)$. We declare $0\in V$ to be the root of $T$ and let
  $\psi:\{1,\ldots,r\}\iso B(T)$ be the identity. Then $(T,0,\psi)$ is
  an ordered tree, and the automorphism $\Pi\iso\Pi$ induced by
  $\sigma\in\Gamma_k$ is compatible with $(T,0,\psi)$, by
  \eqref{hurweq1}. 
\end{exa}

\begin{defn}
  Let $T=(T,v_0,\psi)$ be an ordered tree. A {\em topological
  realization} of $T$ consists of
  \begin{itemize}
  \item
    pairwise distinct points $z_0,\ldots,z_r\in\PP^1_{\CC}$
    (we set $U:=\PP^1_{\CC}-\{z_1,\ldots,z_r\}$),
  \item pairwise disjoint, connected, open subsets $U_v\subset U$, for
    $v\in V(T)$, and
  \item 
    pairwise disjoint simple closed arcs $c_e:[0,1]\to U$, for
    $e\in E(T)$, 
  \end{itemize}
  such that the following holds. The open subsets $U_v$ are the
  connected components of $U-(\cup_e c_e)$. For each edge
  $e=(v_1,v_2)$, the arc $c_e$ lies on the boundary of $U_{v_1}$ and
  $U_{v_2}$, encircling $U_{v_1}$ in clockwise and $U_{v_2}$ in
  counterclockwise direction. For $i=1,\ldots,r$, the set
  $U_{\psi(i)}\cup\{z_i\}$ is homeomorphic to an open disk. Finally,
  $z_0\in U_{v_0}$.  
\end{defn}
  
Fix $v\in V(T)-{v_0}$, let $p_v=(v_0,v_1,\ldots,v_{k-1},v)$ be the
shortest path from $v_0$ to $v$, and set $e:=(v_{k-1},v)$. Set
$z_v:=c_e(0)$ and $z_{v_0}:=z_0$.  Choose a simple path
$a_v:[0,1]\to\bar{U}_v$ leading from $z_{v_{k-1}}$ to $z_v$. Define
\[
    b_v \;:=\; a_{v_1}\cdots a_{v},\quad\text{and}\quad
    \gamma_v \;:=\; b_v\,c_e\,b_v^{-1}.
\]
We regard $\gamma_v$ as an element of $\pi_1(U,z_0)$. To simplify the
notation, set $\gamma_i:=\gamma_{\psi(i)}$, for $i=1,\ldots,r$, and
$\gamma_{v_0}:=1$.

\begin{lem} \label{complem}
  One can choose the paths $a_v$ such that
  \[
       \gamma_v  \;=\; \prod_{v'\in A_v} \gamma_{v'}
  \]
  holds for all $v\in V(T)$ (in particular,
  $\gamma_1\cdots\gamma_r=1$) and that the elements
  $\gamma_1,\ldots,\gamma_r$ induce a presentation
  $\rho:\Pi\iso\pi_1(U,z_0)$.
\end{lem}

A collection $(a_v)$ of paths as in Lemma \ref{complem} will be called
a {\em skeleton} of the topological realization $(z_i,U_v,c_e)$. The
presentation $\rho$ in Lemma \ref{complem} will be said to be induced
by $(a_v)$.

\begin{prop} \label{compprop}
  Let $(z_i,U_v,c_e)$, $(a_v)$ and $\rho$ be as before.
  Let $\tau:U\iso U$ be an orientation preserving diffeomorphism of
  $U$.  Assume that $\tau(z_0)=z_0$ and that for all $v\in V(T)$ 
  the restriction of $\tau$ to $U_v$ is a diffeomorphism $U_v\iso
  U_{v'}$, for some $v'\in V(T)$. Then the induced automorphism of
  $\pi_1(U,z_0)\cong\Pi$ is compatible with $T$.
\end{prop}

\begin{proof}
  It is clear that $\tau$ induces an automorphism $\kappa$ of $T$. We
  will use the notation $v':=\kappa(v)$. Replacing $\tau$ by a
  homotopic diffeomorphism, we may assume that $\tau(c_e)=c_{e'}$, for
  all $e\in E(T)$ (in particular, $\tau(z_v)=z_{v'}$). Fix
  $v\in V(T)-{v_0}$, let $p_v=(v_0,\ldots,v_{k-1},v)$ be the
  shortest path from $v_0$ to $v$ and set $e:=(v_{k-1},v)$. Define
  \[
         \alpha_{v'} \;:=\; 
            b_{v_{k-1}'}\;\tau(a_v)\,a_{v'}^{-1}\;b_{v_{k-1}'}^{-1}
  \]
  and 
  \[
        \beta_{v'} \;:=\; \tau(b_v)\,b_{v'}^{-1} = 
                      \alpha_{v_1'}\ldots \alpha_{v'}
  \]
  Since the closed path $\tau(a_v)\,a_{v'}$ lies entirely in
  $\bar{U}_{v_{k-1}'}$, $\alpha_{v'}$ is an element of the subgroup
  $\Pi_{v_{k-1}'}$. Moreover, we have
  \[
       \tau(\gamma_v) \;=\; \tau(b_v)\,c_{v'}\,\tau(b_v)^{-1}
                      \;=\; \beta_{v'}\,\gamma_{v'}\,\beta_{v'}^{-1}.
  \]
  Therefore, $\tau$ is compatible with $T$ (note that $\chi=1$).
\end{proof}


\subsection{The fundamental group of a degenerating curve} 
\label{adm}

Let $k$ be as in Section \ref{residue} and $S\subset\PP^1_{k^s}$ a
finite, $\Gamma_k$-invariant set and $z_0$ a $k$-rational base point
for $U:=\PP^1_{k^s}-S$. In this section, we associate to $(S,z_0)$ an
ordered tree $T$. The tree $T$ encodes the way in which the points in
$S\cup\{z_0\}$ coalesce on the special fiber. If $k_0$ has
characteristic $0$ then there exists a presentation
$\rho:\Pi\iso\pi_1(U,z_0)$ such that the $\Gamma_k$-action on
$\pi_1(U,z_0)\cong\Pi$ is compatible with the tree $T$, in the sense
of Definition \ref{compdef} above. If $k_0$ is a subfield of $\CC$, we
can actually write down such a presentation.  Moreover, we can compute
the action of the inertia group on $\pi_1(U,z_0)$ explicitly, in terms
of the {\em braid action}.

\vspace{1ex}
The valuation $v$ of $k$ extends uniquely to a valuation on $k^s$,
which we also denote by $v$. Let $R$ be the ring of integers of $k^s$.
We identify the residue field of $R$ with the algebraic closure
$k_0^s$ of $k_0$. Let $(X_R;z_{R,0},\ldots,z_{R,r})$ be the (unique)
model over $R$ of $(\PP^1_{k^s};z_0,\ldots,z_r)$ as a {\em stable
  $(r+1)$-pointed tree of projective lines}, in the sense of
\cite{GHP}. We denote by $\Xb:=X_R\otimes_R k_0^s$ the special fiber.
By definition, $\Xb$ is a tree of projective lines, i.e.\ each
irreducible component of $\Xb$ is non-canonically isomorphic to
$\PP^1$.  For a point $z\in S\cup\{z_0\}$, we denote by $\zb$ its
specialization to $\Xb$.  The points $\zb_0,\ldots,\zb_r$ are pairwise
distinct, smooth points of $\Xb$. Since the model $X_R$ is unique, the
natural $\Gamma_k$-action on $\PP^1_{k^s}$ extends to an action on
$X_R$. This yields an action of $\Gamma_k$ on $\Xb$. In particular,
the inertia group $I\lhd\Gamma_k$ acts $k_0^s$-linearly on $\Xb$.

We define the tree $T$ as follows. The set of vertices $V(T)$ is the
union of $S$ with the set of irreducible components of $\Xb$.  For
every singular point $x$ of $\Xb$, we draw an edge between the two
vertices corresponding to the components meeting in $x$. Moreover, for
each $z\in S$, we draw an edge between the vertex corresponding to $z$
and the vertex corresponding to the component containing $\zb$.  We
declare the component $\Xb_0$ which contains $\zb_0$ as the root of
$T$. It is clear that the elements of $S$ are the leaves of $T$. If
$v\in V(T)$ is not a leaf, we write $\Xb_v$ for the corresponding
component of $\Xb$.  The numbering $\{z_1,\ldots,z_r\}$ of the set $S$
corresponds to a bijection $\psi:\{1,\ldots,r\}\iso S=B(T)$. By
changing this numbering, we may assume that $\psi$ is an order. From
now on, we fix $\psi$ and regard $T=(T,\Xb_0,\psi)$ as an ordered
tree.  We remark that the $\Gamma_k$-action on $\Xb$ induces an action
on $T$. It is clear that this action is determined by the action of
$\Gamma_k$ on $S=B(T)$.

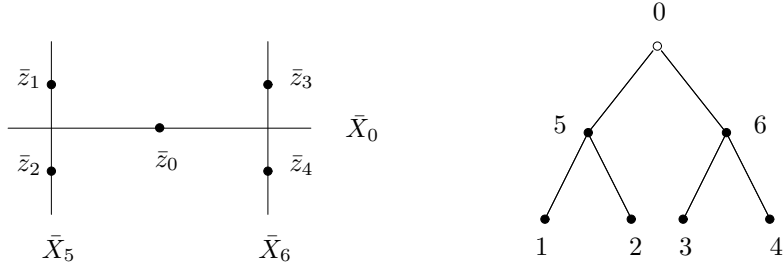
\begin{figure}[bt]
\begin{center}
    \setlength{\unitlength}{0.0005in}

\begin{picture}(8068,3120)(0,-10)
\put(1670,1377){\circle*{90}}
\put(545,1827){\circle*{90}}
\put(545,927){\circle*{90}}
\put(2795,1827){\circle*{90}}
\put(2795,927){\circle*{90}}
\put(6845,2227){\circle{90}}
\put(6125,1327){\circle*{90}}
\put(7565,1327){\circle*{90}}
\put(5675,427){\circle*{90}}
\put(6575,427){\circle*{90}}
\put(7115,427){\circle*{90}}
\put(8015,427){\circle*{90}}
\drawline(95,1377)(3245,1377)
\drawline(545,2277)(545,477)
\drawline(2795,2277)(2795,477)
\drawline(6800,2182)(6125,1327)
\drawline(6800,2182)(6125,1327)
\drawline(6890,2182)(7565,1327)
\drawline(6890,2182)(7565,1327)
\drawline(6125,1327)(5675,427)
\drawline(6125,1327)(5675,427)
\drawline(6125,1327)(6575,427)
\drawline(6125,1327)(6575,427)
\drawline(7565,1327)(7115,427)
\drawline(7565,1327)(7115,427)
\drawline(7565,1327)(8015,427)
\drawline(7565,1327)(8015,427)
\put(185,1827){$\zb_1$}
\put(185,927){$\zb_2$}
\put(3020,1827){$\zb_3$}
\put(3020,927){$\zb_4$}
\put(1620,972){$\zb_0$}
\put(445,27){$\Xb_5$}
\put(2695,27){$\Xb_6$}
\put(3605,1277){$\Xb_0$}
\put(5765,1327){$5$}
\put(7850,1327){$6$}
\put(5575,50){$1$}
\put(6555,50){$2$}
\put(7075,50){$3$}
\put(8015,50){$4$}
\put(6800,2497){$0$}
\end{picture}
\end{center}
\caption{The curve $\Xb$ and the associated tree $T$}
\label{treefig1}
\end{figure}

\begin{exa} \label{admexa1}
  Let $k:=\QQ((t))$, $d\in\ZZ$, $S:=\{\,\pm\sqrt{d}\pm i\sqrt{t}\,\}$
  and $z_0:=0$. The curve $\Xb$ consists of three components $\Xb_0$,
  $\Xb_5$ and $\Xb_6$. Each of these components corresponds to the
  choice of a coordinate $w$ which identifies the function field of
  $\PP^1_{k^s}$ with $k^s(w)$, modulo the action of ${\rm PGL}_2(R)$.
  For the component $\Xb_0$, we can choose the standard coordinate
  $w_0:=z$. For $\Xb_5$, we choose $w_5:=(z+\sqrt{d})/\sqrt{t}$, and
  for $\Xb_6$ we choose $w_6:=(z-\sqrt{d})/\sqrt{t}$.  We let
  $z_1:=-\sqrt{d}+i\sqrt{t}$, $z_2:=-\sqrt{d}-i\sqrt{t}$,
  $z_3:=\sqrt{d}-i\sqrt{t}$ and $z_4:=\sqrt{d}+i\sqrt{t}$. Then
  $z_1,z_2$ reduce to $\zb_1,\zb_2\in\Xb_5$ and $z_3,z_4$ reduce to
  $\zb_3,\zb_4\in\Xb_6$.  See Figure \ref{treefig1}. We remark that
  the standard generator $q_0$ of the inertia group $I$ (which sends
  $t^{1/n}$ to $e^{2\pi i/n}t^{1/n}$) acts as an involution on $\Xb_5$
  and $\Xb_6$.
\end{exa}

The following theorem can be thought of as a rigid-analytic analogue
of Proposition \ref{compprop}. Many similar results can be found in
the literature, see e.g.\ \cite{Nakamura99} or \cite{Pop94}. 

\begin{thm} \label{admthm}
  Assume that $\cha(k_0)=0$. There exists a presentation
  $\rho:\Pi\iso\pi_1(U,z_0)$ such that the $\Gamma_k$-action on $\Pi$
  induced by $\rho$ is compatible with $T$.
\end{thm}

\begin{proof}[Sketch of proof]
  We follow \cite{Saidi97}, with some modifications. The tree $T$ can
  be equipped with a structure $\mathcal{G}$ of a {\em graph of
    groups}. For instance, to each vertex $v\in V(T)$ corresponding to
  a component $\Xb_v$ of $\Xb$ we associate the fundamental group
  $\pi_1(\bar{U}_v)$, where $\bar{U}_v$ is the open subset of $\Xb$
  with all the points $\zb_i$, $i>0$, and all singular points removed.
  To each leaf $\psi(i)\in V(T)$, we associate the group $\ZZh(1)$. We
  obtain a ``canonical'' isomorphism
  \begin{equation} \label{admthmeq}
     \pi_1(U,z_0) \liso \pi_1(T,\mathcal{G}).
  \end{equation} 
  In analogy to Lemma \ref{complem}, one can choose a ``skeleton''
  $(a_v)$ of $(T,\mathcal{G})$ which induces a presentation
  $\Pi_r\iso\pi_1(T,\mathcal{G})$ (actually, to make this precise, one
  has to consider $(T,\mathcal{G})$ as a {\em graphs of groupoids}).
  Let $\rho:\Pi_r\iso\pi_1(U,z_0)$ be the composition of this
  presentation with the isomorphism \eqref{admthmeq}. The rest of the
  proof of Theorem \ref{admthm} is formally the same as the proof of
  Proposition \ref{compprop}.
\end{proof}

\begin{rem} \label{admrem}
  Assume that $\cha(k_0)=p>0$, and let $\pi_1^{\rm adm}(U,z_0)$ be the
  {\em admissible fundamental group}, i.e.\ the inverse limit over the
  finite quotients of $\pi_1(U,z_0)$ corresponding to Galois covers
  with admissible reduction (see \cite{Saidi97}). There exists a
  surjective homomorphism $\rho:\Pi\to\pi_1^{\rm adm}(U,z_0)$ which
  induces an isomorphism on the maximal prime-to-$p$ quotients such
  that the action of $\Gamma_k$ on $\pi_1^{\rm adm}(U,z_0)$ is
  compatible with $\rho$ and $T$, in an obvious sense. In particular,
  Theorem \ref{admthm} remains true if we replace $\Pi$ and
  $\pi_1(U,z_0)$ by their maximal prime-to-$p$ quotients.
\end{rem}

For a more detailed proof of Theorem \ref{admthm} in a special case,
see \cite{boundary}. Under some extra assumptions, we can improve on
Theorem \ref{admthm} and define a concrete presentation $\rho$ with
the claimed properties. This presentation has the advantage that the
inertia action on $\pi_1(U,z_0)$ is known explicitly.

\vspace{1ex}
We assume that $k_0$ is a subfield of the complex numbers, and that
$k:=k_0(t)^{\rm h}$ is the henselization of $k_0(t)$ at the place
$t=0$. For practical applications, this is not a serious restriction.
We will regard $k$ as a subfield of $\CC((t))$ and its algebraic
closure $k^s$ as a subfield of $\CC\{\{t\}\}$. We identify
$\Gamma_{k_0}$ as a subgroup of $\Gamma_k$, via its action on the
coefficients of the Puiseux-expansion of elements of $k$. We let $q_0$
be the ``canonical'' generator of the inertia group $I\lhd\Gamma_k$,
i.e.\ we have $q_0(t^{1/n})=e^{2\pi i/n}t^{1/n}$. For simplicity, we
assume moreover that $\infty\not\in S\cup\{z_0\}$. 

Let $n$ be the smallest positive integer such that all points
$z_0,\ldots,z_r$, regarded as elements of $k^s$, lie in
$\CC((t^{1/n}))$. We let $\tt:=t^{1/n}$ and regard the $z_i$ as
germs of analytic functions of $\tt$. Choose $\epsilon>0$ such that
$z_0,\ldots,z_r$ are meromorphic on the disk
\[
      \Dt \;:=\; \{\; \tt\in\CC \mid\; |\tt|<\epsilon\;\}
\]
and holomorphic on $\Dt^*:=\Dt-\{0\}$. If $\epsilon$ is sufficiently
small, the values $z_i(\tt)$ are pairwise distinct, for all
$\tt\in\Dt^*$. The $R$-curve $X_R$ gives rise to an analytic space
$X_{\Dt}$, together with a map $p:X_{\Dt}\to\Dt$ such that
$p^{-1}(\Dt^*)=\PP^1_{\CC}\times\Dt^*$ and $p^{-1}(0)=\Xb\otimes\CC$.
Let $e=(v_1,v_2)$ be an edge corresponding to a singular point
$x_e\in\Xb$. The complete local ring of $x_e$ on $X_R$ has the form
\[
     \Od_{X_R,x_e} \;=\; R[[\,u_e,v_e\;\mid\; u_ev_e=\tt^{n_e}\;]],
\]
(for some well determined positive integer $n_e$) such that $u_e=0$
(resp.\ $v_e=0$) defines the component $\Xb_{v_1}$ (resp.\ 
$\Xb_{v_2}$) in a neighborhood of $x_e$. We may assume that
$u_e,v_e$ are analytic functions in a neighborhood of $x_e$ on
$X_{\Dt}$. Choosing $\epsilon$ sufficiently small, we may identify the
analytic set
\[
     V_e \;:=\; \{\;(u_e,v_e,\tt)\in\CC^3 \;\mid\;
        u_ev_e=\tt^{n_e},\; |u_e|,|v_e| < \epsilon^{{n_e}/2n}\;\}
\]
with an open neighborhood of $x_e$ on $X_{\Dt}$. We may also assume
that the sets $V_e$ are pairwise disjoint.

Fix a positive real number $t_0$ such that $0<t_0<\epsilon$ and let 
$\zt_i:=z_i(t_0^{1/n})$, for $i=0,\ldots,r$, and
\[
   U_{t_0} \;:=\; \PP^1_{\CC}-\{\zt_1,\ldots,\zt_r\}
           \;\subset\; p^{-1}(t_0^{1/n}).
\]
For each edge $e$ corresponding to a singular point $x_e$, define the
closed arc
\[
    c_e  \;:=\; \left\{
     \begin{array}{ccc}
       [0,1] &  \To    &  U_{t_0}\cap V_e  \\
         s   & \mapsto &  
   (t_0^{n_e/2n}e^{2\pi is},\,t_0^{n_e/2n}e^{-2\pi is},\,t_0^{1/n}) \\
     \end{array} \right..
\]
For the edge $e$ adjacent to the leaf $\psi(i)$, let $c_e:[0,1]\to
U_{t_0}$ be a small closed arc encircling $\zt_i$ in counterclockwise
direction. We may assume that all the arcs $c_e$ are pairwise
disjoint.  For $v\in V(T)-\{v_0\}$, set $e:=(\pre(v),v)\in E(T)$ and
let $U_v$ be the connected component of $U_{t_0}-(\cup_{e'} c_{e'})$
containing the annulus $\{\,(u_e,v_e,\tt)\in V_e \mid
|u_e|>|v_e|\,\}$.

\begin{lem}
   The data $(\zt_i,U_v,c_e)$ is a topological realization of the
   ordered tree $T$.
\end{lem}

Let us choose a skeleton $(a_v)$ for $(\zt_i,U_v,c_e)$ and let
$\rho_{t_0}:\Pi\iso\pi_1(U_{t_0},\zt_0)$ be the induced presentation
(see Lemma \ref{complem}). We define the presentation
$\rho:\Pi\iso\pi_1(U,z_0)$ as the composition of $\rho_{t_0}$ with the
canonical isomorphism
\begin{equation} \label{admeq1}
     \pi_1(U_{t_0},\zt_0) \liso \pi_1(U,z_0).
\end{equation}

A careful modification of the proof of Theorem \ref{admthm} yields
(compare with \cite{Nakamura99}):

\begin{prop} \label{admprop}
  Theorem \ref{admthm} holds with the presentation $\rho$ constructed
  above. 
\end{prop}

Let $0<|t|<\epsilon$ and choose a root $t^{1/n}$. The set
\[
    S_t \;:=\; \{\,z_1(t^{1/n}),\ldots,z_r(t^{1/n})\,\}
        \;\subset\; \PP^1_{\CC}
\]
depends only on $t$, since changing the root $t^{1/n}$ only permutes
the points of $S_t$. Let $\theta:[0,1]\to\CC$ be the closed arc
$\theta(s):=t_0\,e^{2\pi i s}$. The map $s\mapsto S_{\theta(s)}$
corresponds to an element $Q\in \mathcal{B}_r$ of the {\em Artin braid
  group} on $r$ strings. Here we identify $\mathcal{B}_r$ with the
fundamental group of the space of $r$ unordered points in $U_{t_0}$,
with base point $S_{t_0}=\{\zt_1,\ldots,\zt_r\}$. There is a well
known action of $\mathcal{B}_r$ on $\pi_1(U_{t_0},\zt_0)$, see e.g.\ 
\cite{Voelklein}, Section II.10. By the definition of this action, the
braid $Q$ defined above corresponds, via the canonical isomorphism
\eqref{admeq1}, to the inertia generator $q_0$ acting on
$\pi_1(U,z_0)$. This gives us a practical way to compute this action
explicitly, see \cite{Dett1}. 

\begin{exa} \label{admexa2}
  Let $k,S,z_0$ be as in Example \ref{admexa1}. See
  Figure \ref{treefig2} for a picture of the presentation
  $\rho_{t_0}$, for $t_0>0$ sufficiently small and $d>0$. 
  
  It is easily seen that the inertia generator $q_0$ corresponds to
  the braid $Q=Q_1Q_3$ (where $Q_1,Q_2,Q_3$ are the standard
  generators of $\mathcal{B}_r$, see \cite{Voelklein}, Section
  II.10). Therefore, $q_0$ acts on $\Pi$ via
  \begin{equation} \label{admeq3}
     \op{q_0}{\gamma_1}=\gamma_1\gamma_2\gamma_1^{-1},\quad
     \op{q_0}{\gamma_2}=\gamma_1,\quad
     \op{q_0}{\gamma_3}=\gamma_3\gamma_4\gamma_3^{-1},\quad
     \op{q_0}{\gamma_4}=\gamma_3.
  \end{equation}
  Set $\gamma_5:=\gamma_1\gamma_2$, $\gamma_6:=\gamma_3\gamma_4$,
  $\Pi_5:=\gen{\gamma_1,\gamma_2}\subset\Pi$ and
  $\Pi_6:=\gen{\gamma_3,\gamma_4}\subset\Pi$. Let
  $\sigma\in\Gamma_{\QQ}$. By Theorem \ref{admthm} and Proposition
  \ref{admprop}, the action of $\sigma$ on $\Pi\cong\pi_1(U,z_0)$ is
  compatible with the tree $T$. We conclude that
  \begin{equation} \label{admeq4}
    \op{\sigma}{\gamma_i} = \beta_{\sigma(i)}\,
        \gamma_{\sigma(i)}^{\chi(\sigma)}\,\beta_{\sigma(i)}^{-1},
    \qquad\text{with}\quad\beta_i\in
    \begin{cases} \;\Pi_5, & \text{for $i=1,2$} \\
                  \;\Pi_6, & \text{for $i=3,4$.} 
    \end{cases}
  \end{equation}
  Moreover, we have
  \begin{equation} \label{admeq5}
    \op{\sigma}{\gamma_5} =
\gamma_5^{\epsilon(\sigma)\chi(\sigma)},\qquad
    \op{\sigma}{\gamma_6} = \gamma_6^{\epsilon(\sigma)\chi(\sigma)},
  \end{equation}
  where $\epsilon:\Gamma_{\QQ}\to\{\pm 1\}$ is the Kummer character
  $\epsilon(\sigma):=\op{\sigma}{\sqrt{d}}/\sqrt{d}$.
\end{exa}

\begin{figure}[bt]
\begin{center}
    \setlength{\unitlength}{0.0005in}
\begin{picture}(5776,2194)(0,-10)
\put(4688,1092){\circle{2160}}
\put(1088,1092){\circle{2160}}
\put(2888,1092){\circle*{90}}
\put(1088,1542){\circle{540}}
\put(1088,642){\circle{540}}
\put(4688,1542){\circle{540}}
\put(4688,642){\circle{540}}
\put(1088,1542){\circle*{90}}
\put(1088,642){\circle*{90}}
\put(4688,1542){\circle*{90}}
\put(4688,642){\circle*{90}}
\drawline(1013.000,-18.000)(1133.000,12.000)(1013.000,42.000)
\put(1000,-200){$\gamma_5$}
\drawline(4613.000,-18.000)(4733.000,12.000)(4613.000,42.000)
\put(4600,-200){$\gamma_6$}
\drawline(2888,1092)(2168,1092)
\drawline(2888,1092)(3608,1092)(3608,1092)
\drawline(2168,1092)(1358,1542)
\drawline(2168,1092)(1358,642)
\drawline(3608,1092)(4418,1542)
\drawline(3608,1092)(4418,642)
\drawline(1013.000,1242.000)(1133.000,1272.000)(1013.000,1302.000)
\drawline(1013.000,342.000)(1133.000,372.000)(1013.000,402.000)
\drawline(4613.000,1242.000)(4733.000,1272.000)(4613.000,1302.000)
\drawline(4613.000,342.000)(4733.000,372.000)(4613.000,402.000)
\put(2800,1362){$z_0$}
\put(368,1497){$\gamma_1$}
\put(368,642){$\gamma_2$}
\put(5138,1542){$\gamma_4$}
\put(5138,642){$\gamma_3$}
\end{picture}

\end{center}
\caption{}
\label{treefig2}
\end{figure}


\subsection{Computation of the residue}
\label{compute}

In this section we compute the residue $r_v(\omega)$ (and the value
$\omega_0$, if $r_v(\omega)=0$) in a nontrivial example. Let us fix a
square free integer $d$ and let $k:=\QQ((t))$, $S:=\{\pm\sqrt{d}\pm
i\sqrt{t}\}$ and $z_0:=0$, as in Example \ref{admexa1} and Example
\ref{admexa2}. Let $\rho:\Pi\iso\pi_1(U,z_0)$ be the presentation
constructed in the last subsection (see Figure \ref{treefig2}). Then
the action of $\Gamma_k$ on $\pi_1(U,z_0)\cong\Pi$ is subject to
\eqref{admeq3}, \eqref{admeq4} and \eqref{admeq5}.

Let $G:=\At_5$ be the nonsplit central extension of $A_5$ by $C:=\{\pm
1\}$. We denote by $3A$ the (unique) conjugacy class of elements of
order $3$ in $G$ and by $\Ni\inn(3A^4)$ the set of Nielsen classes
$[\g]=[g_1,g_2,g_3,g_4]$ with $g_i\in 3A$, for $i=1,\ldots,4$.

\begin{prop} \label{computeprop}
  Given $\g=(g_1,g_2,g_3,g_4)$ such that $[\g]\in\Ni\inn(3A^4)$, let
  $f:Y\to\PP^1_{k^s}$ be the $G$-cover with branch locus $S$ and
  Hurwitz description $\g$, with respect to the presentation $\rho$.
  We denote by $k'$ the smallest field containing $k$ which is a field
  of moduli for $f$ and by $\omega\in H^2(k',C)$ the
  f.o.d.-obstruction.  Set $g_5:=g_1g_2=(g_3g_4)^{-1}$. Then one of
  the three following cases occurs:
  \begin{enumerate}
  \item 
    Case 1: $g_5=1$. We have $k'=k=\QQ((t))$ and
    $r_v(\omega)=(-d)\in\QQ^*/{\QQ^*}^2$. If $d=-1$ then $r_v(\omega)=0$
    and $\omega_0=(-1,-1)\in\Br_2(\QQ)$.
  \item Case 2: $g_5$ has order $10$. Then
    $k'=\QQ(\sqrt{5})((\tilde{t}))$, with $5v(\tt)=v(t)$,
    $r_v(\omega)=0$ and $\omega_0=(-1,-d)\in\Br_2(\QQ(\sqrt{5}))$.
  \item Case 3: $g_5$ has order $6$. Then $k'=\QQ(\sqrt{-3})((\tt))$,
    with $3v(\tt)=v(t)$, $r_v(\omega)=0$ and
    $\omega_0=(-1,d)=(3,d)\in\Br_2(\QQ(\sqrt{-3}))$.
  \end{enumerate}
\end{prop}

\begin{proof}
  There are exactly $18$ classes $[g_1,g_2,g_3,g_4]$ in
  $\Ni\inn(3A^5)$. Among them, there are two classes with $g_1g_2=1$,
  $10$ classes with $g_5:=g_1g_2$ of order $10$ and $6$ classes with
  $g_5$ of order $6$. See e.g.\ \cite[Prop.~5.8]{BaileyFriedTV} for the
use of a Theorem of Serre \cite{SeLiftAn} showing how to compute the
orders of
products of odd order elements in $\tilde A_n$.
  
  {\bf Case 1:} $g_5=1$. The two Nielsen classes $[\g]\in\Ni\inn(3A^4)$
  such that $g_1g_2=1$ are permuted by an outer automorphism of $G$.
  Therefore, it suffices to consider one of them. So we assume that
  $\g=(g_1,g_1^{-1},g_3,g_3^{-1})$ and that $g_1$ (resp.\ $g_3$) is
  the (unique) lift of $(1\,2\,3)\in A_5$ (resp.\ $(3\,4\,5)\in A_5$)
  to an element of order $3$. By \eqref{admeq3}, we have
  \begin{equation} \label{computeq1}
      \g^{q_0} \;=\; (g_1^{-1},g_1,g_3^{-1},g_3) \;=\;
                     a_0\,\g\,a_0^{-1},
  \end{equation}
  where $a_0\in G$ is a lift of $(1\,2)(4\,5)\in A_5$. Note that $a_0$
  is of order $4$.
  Equation \eqref{admeq4} implies that
\begin{equation} \label{computeq2}
  \g^{\sigma} \;=\; 
  \begin{cases}
    (\,g_1^{\chi_\sigma},\,g_1^{-\chi_\sigma},\,
     g_3^{\chi_\sigma},\,g_3^{-\chi_\sigma}\,), &
          \quad\text{if $\epsilon_\sigma=1$} \\
    (\,g_3^{\chi_\sigma},\,g_3^{-\chi_\sigma},\,
     g_1^{\chi_\sigma},\,g_1^{-\chi_\sigma}\,), &
          \quad\text{if $\epsilon_\sigma=-1$,} 
  \end{cases}
\end{equation}
for all $\sigma\in\Gamma_{\QQ}$ (we have used that the groups
$G_5:=\gen{g_1,g_1^{-1}}$ and $G_6:=\gen{g_3,g_3^{-1}}$ are
cyclic). It follows that $\g^\sigma=b_\sigma\,\g\, b_\sigma^{-1}$,
where
\begin{equation} \label{computeq3}
  b_\sigma \;:=\; 
    a_0^{(1-\chi_\sigma)/2}\;b^{(1-\epsilon_\sigma)/2},
\end{equation}
and $b\in G$ is a lift of $(1\,4)\,(2\,5)\in A_5$. Together with
\eqref{computeq1}, we get that $[\g]^\tau=[\g]$, for all
$\tau\in\Gamma_k$. This shows that $k=\QQ((t))$ is the field of moduli
of the $G$-cover $f$. 

In $\At_5$, we have the equality $b\,a_0\,b^{-1}=a_0^{-1}=-a_0$. Using
this and \eqref{computeq3}, we can compute the $\Gamma_{\QQ}$-action
on $a_0$ (defined by \eqref{actioneq}):
\begin{equation} \label{computeq4}
    \op{\sigma}{a_0} \;:=\; b_\sigma\,a_0\,b_\sigma^{-1}
                     \;=\; a_0^{\epsilon_\sigma}.
\end{equation}
By Proposition \ref{residueprop} (i), the residue $r_v(\omega)\in
H^1(\QQ,\{\pm 1\})$ is represented by the cocycle 
\begin{equation} \label{computeq5}
   c_\sigma \;=\; a_0\,(\op{\sigma}{a_0})^{-\chi_\sigma^{-1}}
            \;=\; a_0^{1-\epsilon_\sigma\chi_\sigma^{-1}}
            \;=\; \epsilon_\sigma\eta_\sigma,
\end{equation}
where $\eta_\sigma:=\op{\sigma}{i}/i$.  Therefore,
$r_v(\omega)=(-d)\in\QQ^*/{\QQ^*}^2=H^1(\QQ,\{\pm 1\})$.  Since we
assumed $d$ to be square free, we have $r_v(\omega)=0$ if and only if
$d=-1$.

If $d=-1$, then $r_v(\omega)=0$ and the value $\omega_0\in
H^2(\QQ,\{\pm 1\})=\Br_2(\QQ)$ of $\omega$ at $v$ is well defined. By
Proposition \ref{residueprop} (ii), $\omega_0$ is the obstruction for
lifting the homomorphism $\phib:\Gamma_{\QQ}\to\Gb=A_5$ given by
\begin{equation} \label{computeq6}
    \phib(\sigma) \;=\; \bar{b}_\sigma \;=\; 
      (1\,2)(4\,5)^{(1-\chi_\sigma)/2}\;
      (1\,4)(2\,5)^{(1-\epsilon_\sigma)/2} 
\end{equation}
to a homomorphism $\phi:\Gamma_{\QQ}\to\At_5$. In other words, it is
the obstruction for lifting the $\ZZ/2\times\ZZ/2$-extension
$\QQ(\sqrt{3},i)/\QQ$ to a quaternion extension. As in the proof of
Proposition \ref{2llprop}, we conclude that
\begin{equation} \label{computeq7}
    \omega_0 \;=\; (-1,3) + (-1,-1) + (3,-1) \;=\; (-1,-1) \;\not=\; 0.
\end{equation}

\vspace{1ex} {\bf Case 2:} $g_5$ has order $10$. Let $O$ be the set of
Nielsen classes $[\g]\in\Ni\inn(3A^4)$ such that $g_5:=g_1g_2$ has
order $10$. The inertia generator $q_0$, which acts on
$\Ni\inn(3A^4)$, stabilizes $O$ and has two orbits $O_1$ and $O_2$ of
length $5$, where a Nielsen class $[\g]\in O$ belongs to $O_1$ (resp.\ 
$O_2$) if $g_5\in 10A$ (resp.\ $10B$). Here $10A$ and $10B$ are the
two conjugacy classes in $\At_5$ of elements of order $10$. Since the
classes $10A$ and $10B$ are conjugate over $\QQ(\sqrt{5})$, it follows
from \eqref{admeq5} that $O_1^\sigma=O_2$ for $\sigma\in\Gamma_{\QQ}$
with $\op{\sigma}{\sqrt{5}}=-\sqrt{5}$.  Therefore, $\Gamma_k$ acts
transitively on $O$, and the fixed field of any class $[\g]\in O$ is
of the form $k':=\QQ(\sqrt{5})((\tt))$, with $5v(\tt)=v(t)$. Hence, all
classes $[\g]$ with $g_5$ of order $10$ are conjugate under the action
of $\Gamma_k$. We may therefore assume that
$\g:=(g_1,g_2,g_2^{-1},g_1^{-1})$, where $g_1$ is a lift of
$(1\,2\,3)\in A_5$ and $g_2$ is a lift of $(1\,4\,5)$. The $G$-cover
$f$ with Hurwitz description $[\g]$ has field of moduli $k'$ as above.
Therefore, there exists an element $a_0\in G$ such that
$\g^{q_0^5}=a_0\,\g\,a_0^{-1}$. In particular,
$a_0\,g_5\,a_0^{-1}=g_5$. Since any subgroup of order $5$ in $A_5$ is
self-centralizing, the image of $a_0$ in $A_5$ has odd order.  From
Corollary \ref{residuecor} we conclude that $r_v(\omega)=0$.

For all $\sigma\in\Gamma_{\QQ(\sqrt{5})}$, there exists $b_\sigma\in
G$ such that $\g^\sigma=b_\sigma\,\g\,b_\sigma^{-1}$. By
\eqref{admeq5}, we have
\begin{equation} \label{computeq8}
   g_5^\sigma \;=\; g_5^{\chi_\sigma\epsilon_\sigma}
                     \;=\; b_\sigma\, g_5\,b_\sigma^{-1}.
\end{equation}
In particular, $b_\sigma$ lies in $N$, the normalizer of
$\gen{g_5}\subset\At_5$. Let $\bar{N}$ be the image of $N$ in $A_5$.
The group $\bar{N}$ is dihedral of order $10$. Let
$\phib(\sigma):=\bar{b}_\sigma$ and
$\psib:\Gamma_{\QQ(\sqrt{5})}\to\gen{\pm 1}$ the composition of
$\phib$ with the sign character $\bar{N}\to\{\pm 1\}$. By
\eqref{computeq8},
$\psib(\sigma)\equiv\chi_\sigma\epsilon_\sigma\pmod{5}$. A computation
shows that $\psib$ corresponds to the quadratic extension
\[
    \QQ\big(\sqrt{5},\sqrt{d(-10+2\sqrt{5})}\;\big)\;/\;\QQ(\sqrt{5}).
\]
To lift $\phib$ to $\At_5$ it suffices to lift $\psib$ to a
homomorphism $\psi:\Gamma_{\QQ(\sqrt{5})}\to\ZZ/4$ (compare with the
proof of Proposition \ref{2llprop}). Using
$10-2\sqrt{5}=(1-\sqrt{5})^2+2^2$, we conclude that
\[
   \omega_0 \;=\; \big(-1,\;d(-10+2\sqrt{5})\,\big) \;=\; (-1,-d).
\]

\vspace{1ex} 
{\bf Case 3:} $g_5$ has order $6$. 
Let $O\subset\Ni\inn(3A^4)$ be the set of classes $[\g]$ with $g_5$ of
order $6$. As in Case 2, $q_0$ acts on $O$ and has two orbits, of
length $3$. Therefore, the fixed field of $[\g]\in O$ is of the form
$k'=k_0((\tt))$, with $3v(\tt)=v(t)$ and $k_0/\QQ$ at most a quadratic
extension. We claim that $k_0=\QQ(\sqrt{-3})$, for all $[\g]\in O$. It
suffices to show the following. Let $\sigma\in\Gamma_{\QQ}$ such that
$\g^\sigma=b_\sigma\,\g\,b_\sigma^{-1}$, for some $b_\sigma\in
G$. Then $\chi(\sigma)\equiv 1\pmod{6}$. 

The subgroups $G_5:=\gen{g_1,g_2}$ and $G_6:=\gen{g_3,g_4}$ of $\At_5$
are of order $24$, isomorphic to $\At_4$ (two $3$-cycles in $A_5$
generate either a cyclic subgroup of order $3$ or a subgroup
isomorphic to $A_4$). Moreover, $G_5\cap G_6$ is cyclic of order $6$,
generated by $g_5$. Now assume that
$\g^\sigma=b_\sigma\,\g\,b_\sigma^{-1}$ and $\epsilon(\sigma)=1$. It
follows from \eqref{admeq4} that $b_\sigma\,G_5\,b_\sigma^{-1}=G_5$
and $b_\sigma\,G_6\,b_\sigma^{-1}=G_6$. It is easy to see that this
implies $b_\sigma\in G_5\cap G_6=\gen{g_5}$. Using \eqref{admeq5}, we
get $g_5^{\chi(\sigma)}=b_\sigma\,g_5\,b_\sigma^{-1}=g_5$, hence
$\chi(\sigma)\equiv 1\pmod{6}$. The argument is similar if
$\epsilon(\sigma)=-1$. Then $b_\sigma\,G_5\,b_\sigma^{-1}=G_6$ and
$b_\sigma\,G_6\,b_\sigma^{-1}=G_5$, and one can show that this implies
$b_\sigma\,g_5\,b_\sigma^{-1}=g_5^{-1}$. Using again \eqref{admeq5},
we obtain $g_5^{-\chi(\sigma)}=b_\sigma\,g_5\,b_\sigma^{-1}=g_5^{-1}$,
hence $\chi(\sigma)\equiv 1\pmod{6}$.

Let $f$ be the $G$-cover with Hurwitz description $[\g]\in O$. We have
shown that the minimal field of moduli of $f$ containing $k$ is of the
form $k'=\QQ(\sqrt{-3})((\tt))$, with $3v(\tt)=v(t)$. Therefore,
$\g^{q_0^3}=a_0\,\g\,a_0$, for some $a_0\in G$. As in Case 2, we
conclude that $a_0\in\gen{g_5}$ and can therefore be assumed to be of
odd order. By Corollary \ref{residuecor}, we have $r_v(\omega)=0$.

Recall that $\g^\sigma=b_\sigma\,\g\,b_\sigma$ implies
$b_\sigma\,g_5\,b_\sigma^{-1}=g_5^{\epsilon(\sigma)}$. Therefore, the
image of $\phib(\sigma):=\bar{b}_\sigma$ is contained in
$\bar{N}\subset A_5$, the image of the normalizer of $g_5$ in $\At_5$.
The group $\bar{N}$ is dihedral of order $6$, and the composition of
$\phib$ with the character $\bar{N}\to\{\pm 1\}$ equals $\epsilon$. As
in Case 2, we conclude that
$\omega_0=(-1,d)=(3,d)\in\Br_2(\QQ(\sqrt{-3}))$.
\end{proof}


\subsection{Specialization to $p$-adic fields}
\label{spec}

In this section, we let $k$ be a field which is complete with respect
to a discrete valuation $v$, and assume that the residue field $k_0$
is perfect of characteristic $p>0$. We denote by $\o_k$ the ring of
integers of $k$ and by $\p$ the maximal ideal of $\o_k$. We set
$K_t:=k((t))$ and denote by $v_t$ the valuation of $K_t$ which has $t$
as a uniformizer. We let $f_t:Y_t\to\PP^1$ be a $G$-cover over
$K_t^s=k^s\{\{t\}\}$ such that $K_t$ is a field of moduli for $f_t$.

Our first goal is to define, for any element $a\in\p-\{0\}$, the 
{\em specialization} $f_a:Y_a\to\PP^1$ of $f_t$ at $t=a$, which should
be
a $G$-cover over $k^s$ with field of moduli $k$. Second, we would like
to compute the residue $r_v(\omega_a)$ of the f.o.d.-obstruction
$\omega_a$ of $f_a$ in terms of $a$ and the residue
$r_{v_t}(\omega_t)$ of the f.o.d.-obstruction $\omega_t$ of $f_t$.
Both these goals are problematic, in general, and we need some extra
assumptions to succeed. 
  
\vspace{1ex} Let $A:=\o_k[[t]]$ and $K:=\Frac(A)$. The ring $A$ is a
complete local domain, regular of dimension $2$ and factorial. The
field $K_t=k((t))$ is the completion of $K$ at the valuation
corresponding to the ideal $(t)\lhd A$. We denote this valuation also
by $v_t$. Moreover, we identify $K^s$ with the algebraic closure of
$K$ inside $K_t^s=k^s\{\{t\}\}$.

\begin{cond} \label{speccond}\ 
  \begin{itemize}
  \item[(a)]
    The branch locus $S_t\subset\PP^1_{K_t^s}$ of the $G$-cover $f_t$
    is defined over $K$, i.e.\ 
    \[
          S_t \;=\; \{ z \;\mid\; F(z)=0\;\},\qquad
          \text{\rm with}\quad F(Z)\in K(Z).
    \]
  (We may choose $F(Z)\in A(Z)$ such that the gcd of the
    coefficients is $1$. We let $\delta(S)\lhd A$ be the ideal
    generated by the discriminant of $F$.)
  \item[(b)]
    We have $\delta(S)=(t^n)$, with $n\geq 0$. 
  \item[(c)]
    The order of the group $G$ is prime to $p=\cha(k_0)$.
  \end{itemize}
\end{cond}

We assume from now on that Condition \ref{speccond} is in force. By
Condition \ref{speccond} (a), the branch locus $S_t\subset
\PP^1_{K_s}$ of the $G$-cover $f_t$ descends to a closed subscheme
$S\subset\PP^1_K$. Let $V:=\Spec A[1/t]$ and $S_V\subset\PP^1_V$ the
closure of $S$ inside $\PP^1_V$. By Condition \ref{speccond} (b),
$S_V\to V$ is \'etale. Therefore, the embedding $S_V\subset\PP^1_V$
corresponds to a morphism $V\to\U_r$, where $\U_r$ is the fine moduli
space for the moduli problem ``$r$ distinct unordered points in
$\PP^1$''. In more concrete terms, we have $\U_r=\PP^r-\delta_r$
($\delta_r$ denotes the discriminant hypersurface), and the
coefficients of the polynomial $F(X)\in A[1/t][X]$ which defines $S_V$
are the projective coordinates for the morphism $V\to\U_r$. Let
$\Hu_r(G)$ be the Hurwitz space classifying $G$-Galois covers over
$\PP^1$ with exactly $r$ branch points. To simplify the notation, we
consider $\U_r$ and $\Hu_r(G)$ as schemes over $\o$. Since the order
of $G$ is prime to the residue characteristic of $\o$ (by Condition
\ref{speccond} (c) ), the natural morphism $\Hu_r(G)\to\U_r$ (which
associates to a $G$-cover its branch locus) is finite \'etale, see
\cite{diss}. The $G$-cover $f_t$ corresponds to a morphism $[f_t]:\Spec
K_t\to\Hu_r(G)$, see Remark \ref{fodrem1} (ii).  

\begin{lem} \label{speclem1}
   The morphism $V\to \U_r$ lifts to a morphism $\phi:V\to\Hu_r(G)$ 
   such that the following diagram commutes.
   \begin{equation} \label{speclem1eq}
   \renewcommand{\arraystretch}{1.5}
   \begin{array}{ccc}
        \Spec K_t      &  \lpfeil{[f_t]} &  \Hu_r(G)  \\
        \big\downarrow & \nearrow & \big\downarrow \\
        V              &  \To     &  \U_r.    \\
   \end{array}
   \end{equation}
\end{lem}

\begin{proof}
  Since $\Hu_r(G)\to\U_r$ is a finite and \'etale morphism, the fiber
  product $V\times_{\U_r}\Hu_r(G)$ will be finite and \'etale over
  $V$. Let $V'\subset V\times_{\U_r}\Hu_r(G)$ be the irreducible
  component which is the image of the morphism $\Spec K_t\to
  V\times_{\U_r}\Hu_r(G)$. Then $V'$ is integral, its fraction field
  $K'$ is a finite extension of $K$ and a subfield of $K_t$, and $V'$
  is the integral closure of $V$ in $K'$. We have to show that $K'=K$.
  
  Let $A'$ be the integral closure of $A$ in $K'$. Clearly, $A'$ is a
  finite $A$-algebra, with fraction field $K'$, and $A'[1/t]$ is
  \'etale over $A[1/t]$. Using Purity of Branch Locus and Abhyankar's
  Lemma, one shows that there exists an embedding $A'\inj
  \o'[[t^{1/n}]]$ of $A$-algebras, for some $n\geq 1$ and some finite
  unramified extension $\o'/\o$. In particular, the valuation $v_t$ of
  $K$ is totally ramified in $K'$. But $K'$ was defined as a subfield
  of $K_t=k((t))$. Therefore, $v_t$ is actually unramified in $K'$.
  Thus, by Purity of Branch locus, $A'/A$ is \'etale. Since $A$ is
  henselian, a finite \'etale extension of $A$ is uniquely determined
  by its residue field extension. It follows that $A'=\o'[[t]]$, for
  some finite unramified extension $\o'/\o$. Using again the embedding
  $A'\subset K_t=k((t))$, we conclude that $A'=A$. This proves the
  lemma.
\end{proof}

Lemma \ref{speclem1} shows in particular that $f_t$ descends to a
$G$-cover $f:Y\to\PP^1_{K^s}$ over $K^s$ such that $K$ is a field of
moduli for $f$, see Remark \ref{fodrem1} (ii).

Let us fix an element $a$ in $\p-\{0\}$. We denote by $v_a$ the
discrete valuation of $K$ corresponding to the ideal $(t-a)\lhd A$
(note that $v_a\not=v_t$). We choose an extension $\tilde{v}_a$ of
$v_a$ to $K^s$. Let $R\subset K^s$ denote the valuation ring
corresponding to $\tilde{v}_a$. By Condition \ref{speccond} (b), the
branch locus $S\subset\PP^1_{K^s}$ of $f$ extends to a subscheme
$S_R\subset\PP^1_R$ such that $S_R$ is \'etale over $\Spec R$. It
follows that $f$ has good reduction at $\tilde{v}_a$. Let
$f_a:Y_a\to\PP^1_{k^s}$ be the reduction of $f$ at $\tilde{v}_a$. We
call $f_a$ the {\em specialization} of $f_t$ at $t=a$. By definition,
$f_a$ is a $G$-cover defined over $k^s$, and $k$ is a field of moduli
for $f_a$. 

The following theorem states that one obtains the f.o.d.-obstruction
for $f_a$ by `specializing' the f.o.d.-obstruction for $f_t$ at
$t=a$. Moreover, this specialization procedure is compatible with
computing the residue. The essential idea behind the theorem is that
the class $\omega_t$ behaves just as one would expect it to behave if
$k$ were a local field of equal characteristic, with field of
coefficients $k_0$, and the class $\omega_t$ arises from pull-back via
the inclusion $k_0((t))\inj k((t))$.

\begin{thm} \label{specthm}
  Assume that Conditions \ref{speccond} (a), (b) and (c) hold. Let
  $\omega_t\in H^2(K_t,C)$ be the f.o.d-obstruction for $f_t$ and
  $\omega_a\in H^2(k,C)$ the f.o.d.-obstruction for $f_a$.  Then
  \begin{enumerate}
  \item
    The residue $r_t(\omega_t)\in H^1(k,C(-1))$ of $\omega_t$ at
    $v_t$ is unramified at $v$, i.e.\ is induced by a class 
    $\bar{r}_t(\omega_t)\in H^1(k_0,C(-1))$. 
  \item
    In $H^1(k_0,C(-1))$, we have the formula
    \[
         r_v(\omega_a) \;=\; v(a)\cdot\bar{r}_t(\omega_t).
    \]
  \item 
    If $r_t(\omega_t)=0$, then $\omega_t$ lies in the submodule
    $H^2(k_0,C)$, and the equality $\omega_a = \omega_t$ holds.
  \end{enumerate}
\end{thm}

\begin{proof}  
  Let $k\nr$ be the maximal unramified extension of $k$ and $\o\nr$
  the ring of integers of $k\nr$. We denote by $\mathcal{S}_A$ the set
  of discrete valuations of $K$ which dominate $A$.  Define
  $L\nr\subset K^s$ as the maximal algebraic extension of $K$ which is
  unramified over each valuation $w\in\mathcal{S}_A$. By Purity of
  Branch Locus, we can identify $\Gal(L\nr/K)$ with $\pi_1(\Spec A)$,
  see \cite{SGA1}.  Since $A$ is a henselian local ring with residue
  field $k_0$, $\pi_1(\Spec A)$ is canonically isomorphic to
  $\Gamma_{k_0}$.  Moreover, $B\nr:=\o\nr[[t]]$ is the integral
  closure of $A$ in $L\nr$.  Similarly, let $L\subset K^s$ be the
  maximal algebraic extension of $K$ which is unramified over each
  valuation $w\in\mathcal{S}_A-\{v_t\}$. We may identify $\Gal(L/K)$
  with $\pi_1(V)$, where $V=\Spec A[1/t]$. Let $B$ be the integral
  closure of $A$ in $L$.  Abhyankar's Lemma together with Purity shows
  that $B=\cup_n\o\nr[[t^{1/n}]]$, where $n$ runs over all integers
  prime to $p$ (compare with the proof of Lemma
  \ref{speclem1}). Therefore, $\Gal(L\nr/L)=\ZZh'(1)$. We obtain the
  following commutative diagram:
  \begin{equation} \label{speceq1}
  {\renewcommand{\arraystretch}{1.5}
  \begin{array}{rcccccl}
    1 \To & \ZZh'(1) & \To & \Gamma_{K_t} & \To & \Gamma_k & \To 1 \\
     & \big\downarrow &  &  \big\downarrow & & \big\downarrow &    \\
    1 \To & \ZZh'(1) & \To & \pi_1(V)     & \To & \Gamma_{k_0} & \To 1.
  \end{array}}
  \end{equation}
  The vertical arrows are induced by the embedding of $K^s$ into
  $K_t^s$. By \cite{SGA1}, we have canonical isomorphism
  $H^i\et(V,C)\cong H^i(\pi_1(V),C)$, for $i\geq 0$. Applying
  \cite[II.A.1, \S 1 and \S2]{SerreCG} to \eqref{speceq1}, we obtain
  the following commutative diagram:
  \begin{equation} \label{speceq2}
  {\renewcommand{\arraystretch}{1.5}
  \begin{array}{rcccccl}
     0\; \to & H^2(k_0,C) & \To & H^2\et(V,C) & 
                    \lpfeil{\bar{r}_t} & H^1(k_0,C(-1)) & \to\; 0 \\
     & \big\downarrow &  &  \big\downarrow & & \big\downarrow &    \\
     0\; \to & H^2(k,C)   & \To & H^2(K_t,C)  & 
                          \lpfeil{r_t} & H^1(k,C(-1)) & \to\; 0.
  \end{array}}
  \end{equation}
  The lower row is nothing else then the residue sequence
  \eqref{reseq} for the field $K_t$. The vertical arrows are the
  restriction homomorphisms from Galois cohomology. 
  
  Let $\phi:V\to\Hu_r(G)$ be the morphism given by Lemma
  \ref{speclem1} and $\omega_V:=\phi^*\tilde{\omega}\in H^2\et(V,C)$
  the pull back of the global f.o.d.-obstruction $\tilde{\omega}\in
  H^2\et(\Hu_r(G),C)$, see Remark \ref{fodrem1} (iii) and \cite{diss}.
  By functoriality and \eqref{speclem1eq}, the morphism
  $H^2\et(V,C)\to H^2(K_t,C)$ in \eqref{speceq2} maps $\omega_V$ 
  to $\omega_t$, the f.o.d.-obstruction for $f_t$. Therefore, Part (i)
  of Theorem \ref{specthm} follows from the commutativity of
  diagram \eqref{speceq2}. 
  
  Let $K_a$ denote the completion of $K$ with respect to the discrete
  valuation $v_a$ and $K_a\nr$ its maximal unramified extension. We
  may identify $K_a$ with $k((t-a))$ and $K_a\nr$ with $k^s((t-a))$.
  Moreover, we may identify $\Gal(K_a\nr/K_a)$ with $\Gamma_k$, the
  absolute Galois group of $k$. Since $v_a$ is unramified in the
  extension $L/K$, there exists an embedding $L\inj
  K_a\nr=k^s((t-a))$. Explicitly, we choose a compatible system
  $a^{1/n}\in k^s$ of $n$th roots of $a\in k$ (where $n$ runs over the
  integers prime to $p$), and define $L\inj K_a\nr$ by
  \begin{equation} \label{speceq3}
     t^{1/n} \;\longmapsto\; a^{1/n}\,(\,1 + \frac{t-a}{a}\,)^{1/n}
             \;=\; a^{1/n}\,(\,1+\frac{t-a}{na}+\cdots\,).
  \end{equation} 
     
  Let $q_0\in\ZZh'(1)$ be a topological generator, corresponding to a
  compatible system $(\zeta_n)$ of $n$th roots of unity, where $p\not|
  n$. Let $\ZZh'(1)\inj \Gamma_k$ be the natural morphism; it sends
  $q_0$ to the automorphism of $k$ which maps $\pi^{1/n}$ to
  $\zeta_n\pi^{1/n}$, where $\pi$ is a uniformizer of $k$.  Similarly,
  let $\ZZh'(1)\inj\pi_1(V)=\Gal(L/K)$ be the morphism that sends
  $q_0$ to the automorphism of $L$ which maps $t^{1/n}$ to
  $\zeta_nt^{1/n}$.  It follows from \eqref{speceq3} that the
  restriction homomorphism $\Gamma_k^t\to\Gal(L/K)=\pi_1(V)$ induces a
  the diagram
  \begin{equation} \label{speceq4}
  {\renewcommand{\arraystretch}{1.5}
  \begin{array}{rcccccl}
    1 \To & \ZZh'(1) & \To & \Gamma_k & \To & \Gamma_{k_0} & \To 1 \\
     & {\scriptstyle\cdot v(a)}\big\downarrow &  &  \big\downarrow & 
                                       & \big\downarrow &    \\
    1 \To & \ZZh'(1)  & \To & \pi_1(V)     & \To & \Gamma_{k_0} & \To 1,
  \end{array}}
  \end{equation}
  where the left vertical arrow sends $q_0$ to $q_0^{v(a)}$.  The
  vertical arrow on the right is clearly the identity.  We apply again
  \cite[II.A.1, \S 1 and \S2]{SerreCG}, this time to \eqref{speceq4},
  and we obtain the following diagram.
  \begin{equation} \label{speceq5}
  {\renewcommand{\arraystretch}{1.5}
  \begin{array}{rcccccl}
     0\; \to & H^2(k_0,C)   & \To & H^2\et(V,C) & 
                    \lpfeil{\bar{r}_t} & H^1(k_0,C(-1)) & \to\; 0 \\
     & \big\downarrow &  &  \big\downarrow & & 
          {\scriptstyle\cdot v(a)}\big\downarrow &    \\
     0\; \to & H^2(k_0,C)   & \To & H^2(k,C)  & 
                          \lpfeil{r_v} & H^1(k_0,C(-1)) & \to\; 0.
  \end{array}}
  \end{equation}
  The upper row is the same as in \eqref{speceq2}, the lower row is
  the residue sequence for the field $k$, see \eqref{reseq}. The left
  and middle vertical arrows are the restriction homomorphism from
  Galois cohomology; in particular, the left vertical arrow is the
  identity. Checking the definition of the residue map (see the proof
  of Proposition \ref{residueprop}), one finds that the vertical arrow
  on the right is multiplication by $v(a)$ and that the diagram
  \eqref{speceq5} commutes. This proves Part (ii) and (iii) of the
  theorem.
\end{proof}

\begin{exa} \label{specexa}
  Let $f:Y\to\PP^1$ be the $\At_5$-covers with branch locus
  $S=\{\pm \sqrt{d}\pm i\sqrt{t}\}$ and Hurwitz description
  $[\g]=[g_1,g_1^{-1},g_3,g_3^{-1}]\in\Ni\inn(3A^4)$ considered in
  Proposition \ref{computeprop}, Case 1. We showed that $\QQ((t))$ is
  a field of moduli for $f$ and that
  $r_t(\omega)=(-d)\in\QQ^*/{\QQ^*}^2$, where $\omega\in
  H^2(\QQ((t)),C)$ is the f.o.d.-obstruction and $r_t$ denotes the
  residue map corresponding to the valuation with parameter $t$.

  Let $p>5$ be a prime and assume that $p\not| d$. Choose an embedding
  of $\QQb$ into $\QQb_p$. We will consider $f$ as a cover over
  $\QQb_p((t))$, with field of moduli $\QQ_p((t))$. The discriminant
  of $S$ is
  \[
     \delta(S) \;=\; (\;2^{12}(d^2+2dt+t^2)d^2t^2\;)
               \;=\; (\,t^2\,) \;\lhd\; \ZZ_p[[t]],
  \]
  and Condition \ref{speccond} holds. We may therefore specialize $f$
  at $t=p$, and we obtain an $\At_5$-cover $f_p:Y_p\to\PP^1$ over
  $\QQb_p$, with branch locus $S_p=\{\pm\sqrt{d}\pm\sqrt{p}\}$ and
  field of moduli $\QQ_p$. Since $p\not|\;d$, the residue
  $r_t(\omega)=(-d)\in H^1(\QQ_p,\pm 1)$ is unramified, as predicted
  by Part (i) of Theorem \ref{specthm}. By Part (ii),
  $r_p(\omega_p)=(-\bar{d})\in H^1(\FF_p,\pm 1)=\FF_p^*/{\FF_p^*}^2$,
  where $\omega_p\in H^2(\QQ_p,\pm 1)=\Br_2(\QQ_p)$ is the
  f.o.d.-obstruction for $f_p$. In terms of the norm residue symbol,
  we have
  \[
       \omega_p \;=\; (\,-d,-p\,)_p \;=\;
       \begin{cases}
           -1 & \text{if $-d\equiv x^2\pmod{p},$} \\
            1 & \text{if $-d\not\equiv x^2\pmod{p}.$}
       \end{cases}
  \]
  It is not hard to extend the computations of Section \ref{compute}
  to show that the f.o.d.-obstruction $\omega_t\in H^2(\QQ((t)),\pm
  1)=\Br_2(\QQ((t))\,)$ is given by the Hilbert symbol $(-d,-t)$. Thus,
  to specialize $\omega_t$ to $\omega_p$, one had to ``plug in'' $t=p$.
\end{exa}


\providecommand{\bysame}{\leavevmode\hbox to3em{\hrulefill}\thinspace}


\end{document}